\documentclass[11pt]{article}

\usepackage{color}
\usepackage{amsmath}
\usepackage{amssymb}
\usepackage{latexsym}

\newtheorem{assumption}{Assumption}
\def\qed{ \ \vrule width.2cm height.2cm depth0cm\smallskip}
\newenvironment{proof}{\noindent {\bf Proof.\/}}{$\qed$\vskip 0.1in}

\newcommand{\la}{\langle}
\newcommand{\ra}{\rangle}

\newcommand{\ba}{\begin{array}}
\newcommand{\ea}{\end{array}}
\newcommand{\be}{\begin{equation}}
\newcommand{\ee}{\end{equation}}
\newcommand{\bea}{\begin{eqnarray}}
\newcommand{\eea}{\end{eqnarray}}
\newcommand{\beaa}{\begin{eqnarray*}}
\newcommand{\eeaa}{\end{eqnarray*}}

\def\neg{\negthinspace}

\def\dbE{\mathbb{E}}
\def\dbF{\mathbb{F}}

\def\dbL{\mathbb{L}}

\def\dbP{\mathbb{P}}
\def\dbR{\mathbb{R}}
\def\dbS{\mathbb{S}}
\def\dbT{\mathbb{T}}
\def\dbQ{\mathbb{Q}}

\def\a{\alpha}
\def\b{\beta}
\def\g{\gamma}
\def\d{\delta}
\def\e{\varepsilon}

\def\l{\lambda}

\def\si{\sigma}

\def\f{\varphi}
\def\th{\theta}
\def\o{\omega}

\def\G{\Gamma}
\def\D{\Delta}
\def\Th{\Theta}

\def\O{\Omega}

\def\uo{\underline\o}
\def\bo{{\boldsymbol{\o}}}
\def\uO{\underline\O}
\def\bO{{\boldsymbol{\O}}}

\def\ab{{\boldsymbol{\a}}}
\def\bd{{\boldsymbol{d}}}

\def\cA{{\cal A}}

\def\cC{{\cal C}}

\def\cF{{\cal F}}

\def\ol{\overline}
\def\ul{\underline}

\def\no{\noindent}

\def\ms{\medskip}
\def\bs{\bigskip}
\def\q{\quad}
\def\qq{\qquad}

\def\pa{\partial}
\def\cd{\cdot}
\def\cds{\cdots}

\def\td{\nabla}

\def\qed{ \hfill \vrule width.25cm height.25cm depth0cm\smallskip}

\newcommand{\basa}{\begin{assumption}}
\newcommand{\easa}{\end{assumption}}

\newcommand{\bas}{\begin{assum}}
\newcommand{\eas}{\end{assum}}

\def\limsup{\mathop{\overline{\rm lim}}}

\def\pa{\partial}
\def\wh{\widehat}

 \def\cd{\cdot}
\def\cds{\cdots}

\def\bnm{{\,|\neg|\neg|\,}}

\def\dis{\displaystyle}

\def\1{{\bf 1}}

\newcommand{\lo}{\stackrel{\leftarrow}{\o}}
\newcommand{\luo}{\stackrel{\leftarrow}{\uo}}
\newcommand{\lbo}{\stackrel{\leftarrow}{\bo}}
\newcommand{\lth}{\stackrel{\leftarrow}{\th}}
\newcommand{\lthpr}{\stackrel{\leftarrow}{\pa_\o \th}}
\newcommand{\lsi}{\stackrel{\leftarrow}{\si}}
\newcommand{\Lg}{\stackrel{\leftarrow}{g}}
\newcommand{\leta}{\stackrel{\leftarrow}{\eta}}

\def\:{\!:\!}
\def\reff#1{{\rm(\ref{#1})}}
\def \proof{{\noindent \bf Proof\quad}}

at 9pt

\begin{document}

\newtheorem{thm}{Theorem}[section]
\newtheorem{lem}[thm]{Lemma}
\newtheorem{cor}[thm]{Corollary}
\newtheorem{prop}[thm]{Proposition}
\newtheorem{rem}[thm]{Remark}
\newtheorem{eg}[thm]{Example}
\newtheorem{defn}[thm]{Definition}
\newtheorem{assum}[thm]{Assumption}

\renewcommand {\theequation}{\arabic{section}.\arabic{equation}}
\def\thesection{\arabic{section}}

\title{\bf  Pathwise It\^{o} Calculus for Rough Paths and Rough PDEs with Path Dependent Coefficients}

\author{Christian  {\sc Keller}\footnote{University of Southern California, Department of Mathematics, kellerch@usc.edu.}   
       \and Jianfeng {\sc Zhang}\footnote{University of Southern California, Department of Mathematics, jianfenz@usc.edu. Research supported in part by NSF grant DMS 1413717.}
\footnote{The authors would like to thank Joscha Diehl, Peter Friz, and Harald Oberhauser for very helpful discussions on the rough path theory and suggestions on the present paper.}}\maketitle

\begin{abstract} 
This paper introduces the path derivatives, in the spirit of Dupire's functional It\^{o} calculus,  for the controlled paths in the rough path theory with possibly non-geometric rough paths.  
The theory allows us to deal with rough integration and rough PDEs in the same manner as standard stochastic  calculus. 
We next study rough PDEs with coefficients depending on the rough path itself, which corresponds to stochastic PDEs with random coefficients. Such coefficients is less regular in the time variable and is not covered in the existing literature.
 The results are useful for studying viscosity solutions of stochastic PDEs.

\end{abstract}

\noindent{\bf Key words:} Rough path, functional It\^{o} calculus, path derivatives, It\^{o}-Ventzell formula, rough differential equations, rough PDEs, stochastic PDEs, characteristics 

\noindent{\bf AMS 2000 subject classifications:} 60H05, 60H10, 60H15,  60G05, 60G17

\vfill\eject

\section{Introduction}
\label{sect-Introduction}
\setcounter{equation}{0}

Firstly initiated by Lyons \cite{Lyons},  the rough path theory has been studied extensively and its applications have been found in many areas, including the recent application on KPZ equations by Hairer \cite{Hairer}. We refer to Lyons \cite{Lyons-ICM}, Friz and Hairer \cite{CF}, Friz and Victoir \cite{FV}, and the reference therein for the general theory and its applications.

On the other hand,  the functional It\^{o} calculus, initiated by Dupire \cite{Dupire} and further developed by Cont and Fournie \cite{CF}, has received very strong attention in recent years. In particular, it has proven to be a very convenient language for viscosity theory of path dependent PDEs, see Ekren, Keller, Touzi and Zhang \cite{EKTZ} and Ekren, Touzi and Zhang \cite{ETZ1, ETZ2}. We also refer to Buckdahn, Ma and Zhang \cite{BMZ},   Cosso and Russo \cite{CR},  Leao, Ohashi and Simas \cite{LOS}, and Oberhauser \cite{Oberhauser} for some recent related works on functional It\^{o} calculus.

The first goal of this paper is to develop the pathwise It\^{o} calculus, in the spirit of Dupire's functional It\^{o} calculus, in the rough path framework with possibly non-geometric rough paths. Based on the bracket process of rough paths, which plays the role of  quadratic variation in semimartingale theory,  we introduce path derivatives for  controlled rough paths of Gubinelli \cite{Gubinelli}. Our first order spatial path derivative is the same as Gubinelli's derivative, and the time derivative is closely related to second order Taylor expansion of the controlled rough paths. This allows us to study the structure of fairly general class of controlled rough paths, and more importantly, to treat the rough integration and rough ODEs/PDEs in the same manner as standard It\^{o} calculus. In particular, as observed by Buckdahn, Ma and Zhang \cite{BMZ} in a Brownian motion setting, we show that the pathwise It\^{o}-Ventzell formula is equivalent to the chain rule of our path derivatives, which is crucial for studying rough PDEs and stochastic PDEs. We shall remark though, while we believe such presentation of path derivatives in rough path framework is new, many related ideas have already been discussed in the literature.   Besides \cite{FH} and the reference therein, we also refer to the recent work Perkowski and Pr\"{o}mel \cite{PP} for some related studies.

We next study the following rough differential equations in the form: 
\bea
\label{RDE0}
d \th_t = g(t, \th_t)  d\bo_t +  f(t, \th_t) d\la\bo\ra_t,
\eea
where $\bo$ is a H\"{o}lder-$\a$ continuous rough path and $\la \bo\ra$ is its bracket process. We remark that we use the Young integration $f(t, \th_t) d\la\bo\ra_t$ rather than Lebesgue integration $f(t, \th_t) dt$ in the drift term above. 
 Our study of above RDE is mainly motivated from the following stochastic differential equations with random coefficients: 
\bea
\label{SDE0}
d X_t = g(t, \o, X_t)  dB_t +  f(t, \o, X_t) dt,
\eea
where $B$ is a Brownian motion in the canonical probability space $(\O,\cF, \dbP)$, $dB$ is It\^{o} integration,  and $g$, $f$ are adapted, namely depend on the history of the path: $\{\o_s\}_{0\le s\le t}$. In the literature, typically the coefficients $g$ and $f$ in \reff{RDE0} do not depend on $t$, or at least is H\"{o}lder-$(1-\a)$ continuous in $t$, see Lejay and Victoir \cite{LV}. However, since a Brownian motion sample path $\o$ is only H\"{o}lder-$({1\over 2}-\e)$ continuous, by setting $\a={1\over 2}-\e$, for \reff{SDE0}  it is not reasonable to assume the mapping $t\mapsto g(\cd, \o, x)$ is H\"{o}lder-$(1-\a)$ continuous as required by  \cite{LV}. Consequently, we are not able to apply the existing results in the rough path literature to study SDE \reff{SDE0} with random coefficients. We shall provide various estimates for rough path integrations, which follow more or less standard arguments, and then establish the wellposedness of RDE \reff{RDE0} under minimum regularity conditions on the coefficients. To be precise, we require only that $g(\cd, x)$, $f(\cd, x)$, and $\pa_\o g(\cd, x)$ are H\"{o}lder-$\b$ continuous for some $\b \in (1-2\a, \a]$,  where $\pa_\o g$ is the spatial path derivative corresponding to Gubnelli's derivative. This can be easily satisfied for the coefficients of \reff{SDE0} when ${1\over 3} < \a < {1\over 2}$.  We note that the recent works  Gubinelli, Tindel and Torrecilla \cite{GTT}, and Lyons and Yang \cite{LY} have also studied the rough integration for more general integrands.  

As a direct consequence of the above wellposedness result of RDE \reff{RDE0}, we obtain the pathwise solution of SDE \reff{SDE0} with random coefficients. Moreover, by restricting the canonical space $\O$ slightly and by using the pathwise stochastic integration, 
we construct the second order process $\uo$ via $\o$ itself. Then the pathwise solution exists for all $\o\in\O$, without the exceptional $\dbP$-null set, and the solution $X(\o)$ is continuous in $\o$ under the rough path topology. 

We would also like to mention that, for linear RDEs, we introduce a decoupling strategy and provide a semi-explicit solution, by using the local solution of certain Riccati type of RDEs. The result seems new even for standard linear SDEs in multidimensional setting.

Finally, we extend the theory to the following rough PDEs with less regular coefficients:
\bea
\label{RPDE0}
d u(t,x) =\big[\si(t,x) \pa_x u + g(t,x, u)\big] d{\bo}_t  + f(t, x, u, \pa_x u, \pa^2_{xx} u)  d\la \bo\ra_t,
\eea
again motivated from pathwise analysis for  stochastic PDEs with random coefficients:
\bea
\label{SPDE0}
d u(t,\o,x) = \big[\si(t,\o, x) \pa_x u + g(t,\o, x, u)\big]  dB_t +  f(t,\o, x, u, \pa_x u, \pa^2_{xx} u) dt.
\eea
As standard in the literature,  see e.g. Kunita \cite{Kunita} for Stochastic PDEs and \cite{FH} for Rough PDEs, the main tool is the  (pathwise) characteristics. We construct the pathwise characteristics via RDEs against a backward rough path. We remark that the backward rough path we construct is also a rough path.  
Our result here is crucial for the study of viscosity solutions of SPDEs in Buckdahn, Ma and Zhang \cite{BMZ2}.

The rest of the paper is organized as follows. In Section 2 we introduce the basics of our pathwise It\^{o} calculus,  in particular  the path derivatives of controlled rough paths. 
In Section 3 we study functions of controlled rough paths and their path derivatives. We shall provide related estimates and prove the chain rule of path derivatives, which is equivalent to the pathwise It\^{o}-Ventzell formula. In Section 4 we study the wellposedness results of rough differential equations. In particular, for linear RDEs we introduce a decoupling strategy which enables us to construct semi-explicit global solution. In Section 5  we apply the RDE results  to SDEs with random coefficients. Finally in Section 6 we extend the results to rough PDEs and stochastic PDEs.

\ms
At below we collect some notations used throughout the paper:

$\bullet$  $T>0$  is a fixed time; and  $\dbT := [0, T]$, $\dbT^2 := \{(s, t): 0\le s < t\le T\}$. 

$\bullet$ $d$ is the fixed dimension for rough paths, and $\dbS^d$ the space of $d\times d$ symmetric matrices.

$\bullet$ $E$ (and $\tilde E$) is a generic Euclid space, and $|E|$ is the dimension of $E$, namely $E = \dbR^{|E|}$.

$\bullet$ By default $E^n$ is viewed as a collum vector. However,  for a function $g: y\in E \to \tilde E$, we take the convention that  the first order derivative  $\pa_y g \in \tilde E^{1\times |E|}$ is viewed as a row vector, and the second order derivative $\pa^2_{yy} g := \pa_y[(\pa_y g)^*]\in \tilde E^{|E|\times |E|}$ is symmetric.  Moreover, for $g: (x,y)\in E_1\times E_2 \to \tilde E$, $\pa_{xy} g := \pa_x[(\pa_y g)^*] \in \tilde E^{|E_2|\times |E_1|}$ and $\pa_{yx} g := \pa_y[(\pa_x g)^*]\in \tilde E^{|E_1|\times E_2}$. 

$\bullet$ $\f_{s, t} := \f_t - \f_s$ for any function $\f:\dbT\to E$ and any $(s, t) \in \dbT^2$.

$\bullet$ For $A\in E^{m\times n}$, $A^*\in E^{n\times m}$ is its transpose.

$\bullet$ For $x\in E^d$ and $y\in \dbR^d$, $ x\cd y \in E$ is their inner product.

$\bullet$ For $A\in E^{m\times n}$ and $\tilde A\in \dbR^{m\times n}$, $ A : \tilde A := \mbox{Trace}(A \tilde A^*) \in E$.

$\bullet$ For $A = [a_{i,j}: 1\le i\le m, 1\le j \le |E|]\in \tilde E^{m\times |E|}$ and $x = [x_{i,j}, 1\le i\le n, 1\le j\le |E|] \in E^n = \dbR^{n\times |E|}$, $A \otimes x\in \tilde E^{m\times n}$ is their convolution whose $(i,j)$-th component is $\sum_{k=1}^{|E|} a_{i,k} x_{j,k}$.

$\bullet$ For $A = [a_{i,j}: 1\le i \le |E_1|, 1\le j \le E_2 ]\in \tilde E^{|E_1|\times |E_2|}$ and $x = [x_{i,j}, 1\le i\le m, 1\le j\le |E_1|] \in E_1 ^m = \dbR^{m\times |E_1|}$, $y = [y_{i,j}, 1\le i\le n, 1\le j\le |E_2|] \in E_2 ^n = \dbR^{n\times |E_2|}$, $A \otimes_2 [x,y]\in \tilde E^{m\times n}$ is their double convolution whose $(i,j)$-th component is $\sum_{k=1}^{|E_1|}\sum_{l=1}^{|E_2|} a_{k,l} x_{i,k} y_{j,l}$.


\section{Rough path integration and path derivatives}
\label{sect-Rough}
\setcounter{equation}{0}
 
In this section we present the basics of rough path theory as well as our pathwise It\^{o} calculus. 
 
\subsection{Rough path and quadratic variation}
 \label{sect-path}
 
 Denote, for a constant $\a>0$,
 \bea
 \label{Oa}
 \left.\ba{c}
 \O_\a (E) :=  \Big\{\o\in C(\dbT, E): \|\o\|_\a <\infty\Big\}, ~\mbox{where}~ \|\o\|_\a := \sup_{(s, t)\in \dbT^2} {|\o_{s,t}|\over |t-s|^\a};\\
 \uO_{\a} (E) :=  \Big\{\uo\in C(\dbT^2, E): \|\uo\|_{\a} <\infty\Big\}, ~\mbox{where}~ \|\uo\|_{\a} := \sup_{(s, t)\in \dbT^2}{|\uo_{s,t}|\over |t-s|^{\a}}.
\ea\right.
\eea
 It is clear that
\bea
\label{oinfty}
\|\o\|_\infty := \sup_{0\le t\le T} |\o_t| \le |\o_0| + T^\a\|\o\|_\a,\q \forall \o\in\O_\a(E).
\eea

From now on, we shall fix two parameters:
 \bea
 \label{ab}
\ab := (\a, \b) &\mbox{where}& \a\in ({1\over 3}, {1\over 2}), \q \b \in (1-2\a, \a] .
 \eea
Our space of rough paths is:
\bea
\label{bOa0}
\bO^0_\a &:=& \Big\{\bo = (\o, \uo)  \in \O_\a (\dbR^d)\times \uO_{2\a}(\dbR^{d\times d}):  \\
&& \uo_{s,t} - \uo_{s, r} - \uo_{r,t} = \o_{s,r}\o^*_{r, t} ~ \forall  0\le s<r<t\le T\Big\}.\nonumber
\eea 
equipped with:
\bea
\label{bOnorm}
\|\bo\|_\a := \|\o\|_\a + \|\uo\|_{2\a}. 
\eea
The requirement in second line of \reff{bOa0} is called Chen's relation. We remark that in general $\|\l \bo\|_\a \neq |\l| \|\bo\|_\a$ for a constant $\l$.

We next introduce the bracket process of $\bo$:
\bea
\label{quadratic}
\la\bo\ra_t := \o_{0,t} (\o_{0, t})^* - \uo_{0, t} - \uo_{0,t}^* \in \dbS^{d}.
\eea

By \reff{bOa0}, one can easily check that
\bea
\label{quadraticproperty}
\la\bo\ra_{s,t} =  \o_{s,t} (\o_{s, t})^*-\uo_{s, t} - \uo_{s,t}^*&\mbox{and thus}& \la\bo\ra \in \O_{2\a}(\dbS^d).
\eea

\begin{rem}
\label{rem-quadratic}
{\rm (i) Clearly $\la \bo\ra=0$ if and only if $\bo$ is a geometric rough path. This process  is intrinsic for non-geometric rough paths, and makes our  study  much more convenient. 

(ii) The process $\la \bo\ra$ is called the bracket process, denoted as $[\bo]$, of the so called reduced rough path in \cite{FH}.  As we will see later, this process plays essentially the same role as the quadratic variation process in semimartingale theory. However, we shall note that a typical rough path may not have finite  quadratic variation.
\qed}
\end{rem}

The following result is straightforward and its proof is omitted.

\begin{lem}
\label{lem-quadratic}
For any $\bo, \tilde\bo \in \bO_\a^0$, we have
\bea
\label{Qonorm}
\|\la\bo\ra\|_{2\a} \le \|\bo\|_\a[2+\|\bo\|_\a];\q
\|\la\bo\ra - \la \tilde \bo\ra\|_{2\a} \le   [\|\o\|_\a+\|\tilde\o\|_\a+2] \|\bo- \tilde\bo\|_\a.
\eea
\end{lem}

\subsection{Rough path integration}
 \label{sect-integration}
 
To study rough path integration against $\bo$, we first  introduce the controlled rough paths of Gubinelli \cite{Gubinelli}, which can be viewed as $C^1$-regularity of the paths against the rough path.  
 
 \begin{defn}
 \label{defn-Gubinelli}
 For each $\o \in \O_\a(\dbR^d)$,  the space $\cC^1_{\o,\ab}(E)$  consists of  $E$-valued controlled rough paths $\th\in \O_\b(E)$ such that there exists $\pa_\o \th \in \O_\b(E^{1\times d})$ satisfying:
 \beaa
 R^{\o,\th}\in \uO_{\a+\b}(E)
&\mbox{where}&   R^{\o, \th}_{s,t} :=  \th_{s,t} - \pa_\o\th_s  \o_{s,t}, \forall (s,t)\in \dbT^2.
\eeaa
\end{defn}
We note that for notational simplicity we take the convention that $\pa_\o \th$ is a row vector.

\begin{rem}
\label{rem-Gubinelli}
{\rm (i) The $\pa_\o \th$ depends on $\o$, but not  on $\uo$. 

(ii) In general $\pa_\o \th$ is not unique. However, when $\bo$ is truly rough, namely $\bo \in \bO_\a$ as defined in \reff{bOa} below, $\pa_\o \th$ is unique. See \cite{FH} Proposition 6.4. For the ease of presentation, in this paper we shall assume $\bo \in \bO_\a$. However, most of our results still hold true when $\bo \in \bO^0_\a$, provided that we specify a version of $\pa_\o \th$.

(iii) $\pa_\o \th$ is called the Gubinelli derivative in the rough path literature. As we will see in Section \ref{sect-SDE}, when $\o$ is a sample path of Brownian motion, it coincides with the path derivative introduced in \cite{BMZ}. So in this paper we also call it path derivative.
\qed}
\end{rem}
For the ease of presentation, from now on we restrict to $\bo\in \bO_\a$ so that $\pa_\o\th$ is unique:
\bea
\label{bOa}
\bO_\a &:= \Big\{\bo \in \bO^0_\a:&  \mbox{ there exists a dense subset $A\subset [0, T)$ such that}\\
&&  \limsup_{t\downarrow s} {|v\cd \o_{s,t}|\over (t-s)^{2\a}} = \infty~~\mbox{for all $s\in A$ and $v\in \dbR^d \backslash \{0\}$}\Big\}.\nonumber
\eea 
For $\bo\in \bO_\a$, we equip the space $\cC^1_{\o,\ab}(E)$ with the semi-norms: 
\bea
\label{Canorm}
\left.\ba{c}
\|\th\|_{\o,\ab} :=\|\pa_\o\th\|_\b + \|R^{\o, \th}\|_{\a+\b},\q d^{\o,\tilde\o}_{\ab}(\th, \tilde\th) := \|\pa_\o\th-\pa_{\tilde\o}\tilde\th\|_\b + \|R^{\o, \th} - R^{\tilde\o, \tilde\th}\|_{\a+\b},\\
\bnm \th\bnm_{\o,\ab} := \|\th\|_{\o,\ab}  + |\pa_\o\th_0|,\q \bd^{\o,\tilde\o}_{\ab}(\th, \tilde\th) := d^{\o,\tilde\o}_{\ab}(\th, \tilde\th) + |\pa_\o \th_0 - \pa_{\tilde\o}\tilde\th_0|.
\ea\right.
\eea
In particular, we note that
\bea
\label{Canorm2}
 d^{\o}_{\ab}(\th, \tilde\th) := d^{\o,\o}_{\ab}(\th, \tilde\th) = \| \th- \tilde\th\|_{\o, \ab},~~  \bd^{\o}_{\ab}(\th, \tilde\th) := \bd^{\o,\o}_{\ab}(\th, \tilde\th) = \bnm  \th-\tilde\th \bnm_{\o, \ab}.
 \eea
By  \reff{oinfty} one can easily check that
\bea
\label{tha}
\left.\ba{c}
\O_{\a+\b}(E) \subset \cC^1_{\o,\ab}(E),\q \mbox{with}~ \pa_\o\th=0 ~\mbox{and}~ \|\th\|_{\o,\ab} = \|\th\|_{\a+\b},  ~\forall \th\in \O_{\a+\b};\\
\cC^1_{\o,\ab}(E) \subset \O_\a(E),\q \mbox{with}~  \|\th\|_\a \le |\pa_\o \th_0|\|\o\|_\a + T^\b[1+\|\o\|_\a] \|\th\|_{\o,\ab}~\forall \th\in \cC^1_{\o,\ab}(E).
\ea\right.
\eea

We are now ready to define the rough path integration. For each $\bo\in \bO_\a$, $\th \in \cC^1_{\o,\ab}(E^d)$, and each partition $\pi: 0=t_0<\cds<t_n=T$, denote
\bea
\label{Thpi}
\Th^\pi_t :=  \sum_{i=0}^{n-1}\Big[\th_{t_i\wedge t} \cd \o_{t_i\wedge t, t_{i+1}\wedge t}  + \pa_\o \th_{t_i\wedge t} : \uo_{t_i\wedge t, t_{i+1}\wedge t}\Big].
\eea
Here, for $\th = [\th_1,\cds, \th_d]^*$, we take the convention that $\pa_\o \th\in E^{d\times d}$ with $i$-th row $\pa_\o\th_i$. 
Following Gubinelli \cite{Gubinelli}, we may define the rough integral as the unique limit of $\Th^\pi$:
\begin{lem}
\label{lem-RI} For each $\bo\in \bO_\a$, $\th \in \cC^1_{\o,\ab}(E^d)$, the rough integral
\bea
\label{RI}
\int_0^t \th_s \cd d\bo_s :=\Th_t:=  \lim_{|\pi|\to 0} \Th^\pi_t \in E
\eea
exists, and is independent of the choice of $\pi$. Moreover, $\Th\in \cC^1_{\o,\ab}(E)$ with $\pa_\o \Th = \th^*$ and 
\bea
\label{RIest}
\left.\ba{c}
\dis \Big|\Th_{s,t} - \th_s \cd \o_{s,t} - \pa_\o \th_s : \uo_{s,t}\Big| \le  C_{\ab}\|\bo\|_\a \|\th\|_{\o, \ab} |t-s|^{2\a+\b};\\
\dis \|\Th\|_{\o,\ab} \le T^{\a-\b}\|\bo\|_{\a}  |\pa_\o \th_0|+ C_{\ab} T^\a[1+\|\bo\|_\a] \|\th\|_{\o, \ab},
\ea\right.
\eea
where the constant $C_{\ab} $ depends only on $\ab$ and the dimensions $|E|$ and $d$. 
\end{lem}
\proof This result follows the same arguments in \cite{FH} Theorem 4.10, except that the second line of \reff{RIest} appears slightly differently. To see that, by the first estimate we have
\beaa
 \|R^{\o,\th}\|_{\a+\b} \le  \|\pa_\o \th\|_\infty \|\bo\|_{\a} T^{\a-\b} + CT^\a\|\bo\|_\a \|\th\|_{\o, \ab}.
\eeaa
Plug the first inequality of \reff{tha} into above and then use the second inequality of \reff{tha}, we obtain the second estimate of \reff{RIest} immediately. 
\qed

Moreover, we have the following stability result in terms of the rough integral, which improves \cite{FH} Theorem 4.16 slightly.
\begin{lem}
\label{lem-DRI}
Let $(\bo, \th, \Th)$ be as in Lemma \ref{lem-RI} and consider $(\tilde\bo, \tilde\th, \tilde \Th)$ similarly. Denote 
\beaa
M :=  \|\th\|_{\o,\ab} + \|\tilde\th\|_{\tilde\o, \ab} + \|\bo\|_\a+\|\tilde\bo\|_\a, &\mbox{and}& \D \f := \tilde\f-\f,~\mbox{for}~ \f = \bo, \th, \Th.
\eeaa 
Then, there exists a constant $C_{\ab,M}$, depending on $\ab, M$, and $|E|$, $d$, such that
\bea
\label{DRIest}
d^{\o,\tilde\o}_{\ab}(\Th,  \tilde\Th) \le T^{\a-\b} \Big[ | \pa_\o \tilde \th_0|  \|\D\bo\|_\a + \|\bo\|_{\a}  |\D\pa_\o \th_0| \Big]+  C_{\ab,M} T^\a\Big[  \|\D\bo\|_\a +  d^{\o,\tilde\o}_{\ab}( \th,  \tilde\th)\Big].\nonumber
\eea
\end{lem}
\proof First, similar to the first estimate in \reff{RIest}, or following the same arguments as in  \cite{FH} Theorem 4.16, we have
\beaa
\Big|[R^{\tilde\o, \tilde\Th}_{s,t} -\pa_\o \tilde \th_s : \tilde \uo_{s,t}]-[R^{\o,\Th}_{s,t}  - \pa_\o \th_s : \uo_{s,t}] \Big|
\le CT^\a\Big[  \|\D\bo\|_\a +  d^{\o,\tilde\o}_{\ab}( \th,  \tilde\th)\Big](t-s)^{\a+\b}.
\eeaa
Note that, by \reff{oinfty},
\beaa
&& |\pa_\o \tilde \th_s : \tilde \uo_{s,t}-\pa_\o \th_s : \uo_{s,t} | \le \Big[\|\D\pa_\o \th\|_\infty \|\uo\|_{2\a} + \|\pa_\o \tilde\th\|_\infty \|\D\uo\|_{2\a}\Big] (t-s)^{2\a}\\
&&\le  \Big[[|\D\pa_\o \th_0| \|\uo\|_{2\a} + |\pa_\o \tilde\th_0| \|\D\uo\|_{2\a} ]+ C T^\b [\|\D\pa_\o \th\|_\b  + \|\D\uo\|_{2\a}]\Big] (t-s)^{2\a}.
\eeaa
Then we obtain the desired estimate for $\| R^{\tilde\o, \tilde \Th}-R^{\o,\Th}\|_{\a+\b}$ immediately. Moreover,
\beaa
|\D\pa_\o \Th_{s,t}| &=& |\D\th_{s,t}| = \Big|[\pa_\o \tilde\th_s \tilde\o_{s,t} + R^{\tilde \o,\tilde \th}_{s,t}]-[\pa_\o \th_s \o_{s,t} + R^{\o,\th}_{s,t}]\Big|\\
&\le&\Big[\|\D\pa_\o \th\|_\infty \|\o\|_{\a} + \|\pa_\o \tilde\th\|_\infty \|\D\o\|_{\a} + T^\b\|R^{\tilde\o, \tilde\th}-R^{\o,\th}\|_{\a+\b}\Big] (t-s)^{\a}
\eeaa
By \reff{oinfty} again we obtain the  desired estimate for $\|\D\pa_\o \Th\|_\b$, completing the proof.
\qed

We conclude this subsection with the Young's integration against $\la \bo\ra$. Since $\la\bo\ra\in \O_{2\a}(\dbS^d)$, by \reff{ab}  the Young's integral $\th_t : d\la\bo\ra_t$ is well defined  for all $\th \in \O_\b(E^{d\times d})$. We collect below some results concerning  this integration. Since the proofs are standard and are  much easier than Lemmas \ref{lem-RI} and \ref{lem-DRI}, we thus omit them.

\begin{lem}
\label{lem-Young} (i)  Let $\bo\in \bO_\a$, $\th \in \O_\b(E^{d\times d})$,   $\Th_t := \int_0^t \th_s :d\la\bo\ra_s$. Then  $\Th\in \O_{\a+\b}(E)$ and 
\bea
\label{Young}
\left.\ba{c}
|\Th_{s,t} - \th_s: \la \bo\ra_{s,t}| \le C\|\th\|_\b \|\la\bo\ra\|_{2\a} (t-s)^{2\a+\b}, \\
\|\Th\|_{\a+\b} \le \Big[T^{\a-\b}|\th_0| + CT^\a \|\th\|_\b \Big]\|\la\bo\ra\|_{2\a}.
\ea\right.
\eea
(ii) Let $(\tilde \bo, \tilde\th, \tilde \Th)$ satisfy the same properties. Then, denoting $\D \f := \f-\tilde \f$ for $\f = \bo, \th, \Th$,
\bea
\label{DYoung}
\|\D\Th\|_{\a+\b} \le T^{\a-\b}\|\la\bo\ra\|_{2\a}|\D\th_0| + CT^\a\Big[\|\la\bo\ra\|_{2\a} \|\D \th\|_\b + \|\tilde\th\|_\b \|\la\bo\ra - \la \tilde \bo\ra\|_{2\a}\Big].
\eea
\end{lem}

\subsection{Path derivatives}
We next introduce further path derivatives of $\th$. Our following definition is motivated from the path derivatives introduced in Ekren, Touzi and Zhang \cite{ETZ1} and Buckdahn, Ma and Zhang \cite{BMZ}, which in turn were motivated by the functional It\^{o} calculus of Dupire \cite{Dupire}.

 \begin{defn}
 \label{defn-C2}
 For each $\bo \in \bO_\a$,  the space $\cC^{2}_{\bo,\ab}(E)$  consists of  $E$-valued controlled rough paths $\th\in \cC^{1}_{\o,\ab}(E)$ such that $\pa_\o \th \in \cC^{1}_{\o,\ab}(E^{1\times d})$ and there exists symmetric $D^\bo_t \th \in \O_\b(E^{d\times d})$  satisfying the following pathwise It\^{o} formula:
 \bea
 \label{Ito}
 d \th_t = \pa_\o \th_t d\bo_t + [D^\bo_t \th_t + {1\over 2} \pa^2_{\o\o} \th_t ] : d\la\bo\ra_t, ~\mbox{where}~ \pa^2_{\o\o} \th_t := \pa_\o [(\pa_\o \th_t)^*]\in E^{d\times d}
\eea
\end{defn}

\begin{rem}
\label{rem-Dt}
{\rm  (i) In general $D^\bo_t \th$ may not be unique. Similar to \reff{bOa}, one can easily check that $D^\bo_t\th$ is unique if $\bo$ is restricted to the following $\wh \bO_\a$:
 \bea
\label{hatbOa}
\wh\bO_\a &:= \Big\{\bo \in \bO_\a:&  \mbox{ there exists a dense subset $A\subset [0, T)$ such that}\\
&&  \limsup_{t\downarrow s} {|v : \la\bo\ra_{s,t}|\over (t-s)^{2\a+\b}} = \infty~~\mbox{for all $s\in A$ and $v\in \dbS^d \backslash \{0\}$}\Big\}.\nonumber
\eea 

(ii) However, $\la \bo\ra$ is more regular than $\bo$, and thus \reff{hatbOa} is much more difficult to satisfy than \reff{bOa}. For example, if $\bo$ is a sample path of Brownian motion with It\^{o} integration, then $\la \bo\ra_t = tI_d$ as we will see in Section \ref{sect-SDE} below. Consequently, by considering $v   \in \dbS^d \backslash \{0\}$ with Trace$(v) =0$, we see that $\wh\bO_\a = \emptyset$. 


 (iii) In many cases in this paper, $\th$ already takes the form $d\th_t = a_t \cd d\bo_t + b_t : d\la\bo\ra_t$,  then clearly $ \pa_\o \th = a^*$ and we shall always set,  thanks to the symmetry of $\la\bo\ra$,
 \bea
 \label{Dbot}
D^\bo_t \th := {1\over 2}\Big[ (b- {1\over 2}\pa_\o a) + (b-{1\over 2}\pa_\o a)^*\Big].
 \eea

(iv) In the case that $\la \bo\ra_t = t$, we will actually define $\pa^\bo_t \th := $Trace$(D^\bo_t\th)$. Then we see that $\pa^\bo_t\th$ is unique. 
\qed}
\end{rem}

\begin{rem}
\label{rem-C2}
{\rm 
(i) In general $\pa_{\o^i}$ and $\pa_{\o^j}$ do not commute, and  $D^\bo_t$ and $\pa_\o$ are also not commutative.  In particular, $\pa^2_{\o\o} \th$ is not symmetric. However, since $\la\bo\ra$ is symmetric, we see that \reff{Ito} is equivalent to
 \bea
 \label{Ito2}
 d \th_t = \pa_\o \th_t d\bo_t + \Big[D^\bo_t \th_t + {1\over 4} [\pa^2_{\o\o} \th_t + (\pa^2_{\o\o}\th_t)^*]\Big] : d\la\bo\ra_t.
 \eea

(ii) One can easily check that the pathwise It\^{o} formulae \reff{Ito} and \reff{Ito2} are equivalent to the following pathwise Taylor expansion:
 \bea
 \label{Taylor1}
 \th_{s,t} = \pa_\o \th_s \o_{s,t} + {1\over 2} \pa^2_{\o\o} \th_s  : [\o_{s,t} \o^*_{s,t} + \uo_{s,t} - \uo_{s,t}^*] + D^\bo_t \th_s : \la\bo\ra_{s,t} + O((t-s)^{2\a+\b}).
 \eea
 In the case that $\pa^2_{\o\o} \th$ is symmetric, which is always the case when $d=1$, \reff{Taylor1} becomes
 \bea
 \label{Taylor2}
 \th_{s,t} = \pa_\o \th_s \o_{s,t} + {1\over 2} \pa^2_{\o\o} \th_s  : [\o_{s,t} \o^*_{s,t}] + D^\bo_t \th_s : \la\bo\ra_{s,t} + O((t-s)^{2\a+\b}).
 \eea
 We refer to \cite{BMZ} for related works in Brownian motion setting.
\qed}
\end{rem}

\subsection{Backward rough integration}
 In this subsection we introduce the backward rough path, which is also a rough path and will play an important role in constructing the pathwise characteristics in Section \ref{sect-RPDE} below.  Let $\bo \in \bO_\a$ and $\th\in \cC^1_{\o,\ab}(E^d)$. For any $t_0 \in [0, T]$ and $0\le s \le t\le t_0$, define
\bea
\label{back}
\left.\ba{c}
\lo^{t_0}_t := \o_{t_0}-\o_{t_0-t},\q \luo^{t_0}_{s,t} :=   \o_{t_0-t, t_0-s} \o^*_{t_0-t, t_0-s} - \uo_{t_0-t, t_0-s},\q \lbo^{t_0} := (\lo^{t_0},  \luo^{t_0});\\
 \lth^{t_0}_t := \th_{t_0-t},\q (\lthpr)^{t_0}_t := -\pa_\o \th_{t_0-t}.
\ea\right.
\eea
 By restricting the processes on $[0, t_0]$ in obvious sense, we have
\begin{lem}
\label{lem-Back}
Let $\bo \in \bO_\a$ and $\th\in \cC^1_{\o,\ab}(E^d)$. Then $\lbo^{t_0} \in  \bO^0_\a$, $\lth^{t_0} \in \cC^1_{\lo^{t_0},\a}(E^d)$ with 
\bea
\label{BI}
\pa_{\lbo^{t_0}} \lth^{t_0} = (\lthpr)^{t_0} &\mbox{and}& \int_{t_0-t}^{t_0-s} \lth^{t_0}_r \cd d\lbo^{t_0}_r = \int_s^{t} \th_r\cd d \bo_r,\q 0\le s<t\le t_0.
\eea
\end{lem}
\proof In this proof we omit the superscript $^{t_0}$ and denote $t':= t_0-t$, $s':= t_0-s$, $r':= t_0-r$,  $\d := t-s$. First, one can easily check that 
\beaa
\lo_{s,t} = \o_{t', s'},\q \luo_{s,t}-\luo_{s,r}-\luo_{r,t} = \o_{r', s'} \o^*_{t', s'} = \lo_{s,r}\lo_{r,t}.
\eeaa
This implies that $\lbo \in \bO^0_\a$. Next,
\beaa
\lth_{s,t} = -\th_{t', s'} = - \pa_\o \th_{t'}\o_{t', s'} - R^{\o, \th}_{t', s'} = \lthpr_s \lo_{s,t}  + \pa_\o \th_{t', s'} \o_{t', s'} - R^{\o, \th}_{t', s'}.
\eeaa
Then clearly $ \lthpr$ is a Gubinelli derivative of $\lth$ with respect to $\lo$. Finally, the second equality of \reff{BI} is exactly the same as \cite{FH} Proposition 5.10.
\qed

We remark that $\lbo^{t_0}$ may not be in $\bO_\a$, and then  $\pa_{\lbo^{t_0}} \lth^{t_0}$ is not unique. See Remark \ref{rem-Gubinelli} (ii). In this case we shall always choose $(\lthpr)^{t_0}$ as its path derivative.

\section{Functions of controlled paths} 
\label{sect-function}
\setcounter{equation}{0}

In this section we study functions $\f: \dbT\times \tilde E \to  E$ and its related path derivatives. Similar to \reff{Ito}, we shall take the notational convention that
\bea
\label{yo}
\pa_{yy} \f := \pa_y [(\pa_y \f)^*],\q \pa_{y\o} \f := \pa_y [(\pa_\o\f)^*],\q \pa_{\o y}\f := \pa_\o [(\pa_y \f)^*].
\eea


\begin{defn}
\label{defn-f}

(i) For $k\ge 0$, let  $\cC^{k}_{ loc}(\tilde E, E)$ be the set of mappings $g: \dbT\times \tilde E \to E$ such that $g$ is $k$-th differentiable in $y$.  Moreover, let  $\cC^{k}(\tilde E, E)\subset \cC^{k}_{loc}(\tilde E, E)$ be such that 
\bea
\label{gk}
\|g\|_{k} := \sum_{i=0}^k \sup_{y \in  \tilde E} \|\pa_y^{(i)} g(\cd,y)\|_\infty <\infty.
\eea

(ii) For $k\ge 0$,  let  $\cC^{k}_{\b, loc}(\tilde E, E)\subset \cC^{k}_{loc}(\tilde E, E)$ be such that, for $i=0,\cds, k$, $\pa^{(i)}_y g$ is H\"{o}lder-$\b$ continuous in $t$, and the mapping $y\mapsto \pa^{(i)}_y g(\cd, y)$ is continuous under $\|\cd\|_\b$. Moreover,  let  $\cC^{k}_{\b}(\tilde E, E)\subset \cC^{k}_{\b, loc}(\tilde E, E)$ be such that
\bea
\label{gkb}
\|g\|_{k,\b} := \sum_{i=0}^k \sup_{y\in \tilde E} \|\pa_y^{(i)} g(\cd,y)\|_\b <\infty.
\eea

(iii) Let  $\cC^{1,2}_{\o,\ab, loc}(\tilde E, E)\subset \cC^{2}_{loc}(\tilde E, E)$ be such that  $g(\cd,y)\in \cC^1_{\o,\a}(E)$,  $\pa_y g(\cd, y) \in \cC^{1}_{\o,\ab}(E^{1\times |\tilde E|})$, for each $y\in \tilde E$, the mappings $y\mapsto g(\cd, y)$ and $y\mapsto \pa_y g(\cd, y)$ are continuous under $\bnm\cd\bnm_{\o,\ab}$, and $\pa_\o g \in \cC^1_{\b, loc}(\tilde E, E^{1\times d})$. Moreover, let  $\cC^{1,2}_{\o,\ab}(\tilde E, E)\subset \cC^{1,2}_{\o,\ab, loc}(\tilde E, E)$ be such that 
\bea
\label{gC2}
\|g\|_{2,\bo,\ab} := \|g\|_2 + \|\pa_\o g\|_1 + \sup_{y\in \tilde E} [ \|g(\cd, y)\|_{\bo,\ab} +    \|\pa_y g(\cd, y)\|_{\bo,\ab}] <\infty.
\eea

(iv) Let  $\cC^{2,3}_{\o,\ab, loc}(\tilde E, E)\subset \cC^{1,2}_{\o,\ab, loc}(\tilde E, E)$ be such that $\pa_\o g \in \cC^{1,2}_{\o,\ab, loc}(\tilde E, E^{1\times d})$, $\pa_y g \in \cC^{1,2}_{\o,\ab, loc}(\tilde E, E^{1\times |\tilde E|})$, $g(\cd, y)\in \cC^2_{\o,\ab}(E)$ for every $y\in \tilde E$ and there exists $D^\o_t g \in \cC^1_{\b,loc}(\tilde E, E^{d\times d})$. Moreover, let  $\cC^{2,3}_{\o,\ab}(\tilde E, E)\subset \cC^{2,3}_{\o,\ab, loc}(\tilde E, E)$ be such that 
\bea
\label{gC3}
\|g\|_{3,\bo,\ab} :=\|g\|_{2,\bo,\ab} + \|\pa_\o g\|_{2,\bo,\ab}+ \|\pa_y g\|_{2,\bo,\ab} <\infty.
\eea

(v) Let  $\cC^{3,3}_{\o,\ab, loc}(\tilde E, E)\subset \cC^{2,3}_{\o,\ab, loc}(\tilde E, E)$ be such that $\pa_\o g \in \cC^{2,3}_{\o,\ab, loc}(\tilde E, E^{1\times d})$.

(vi) For $\bo, \tilde \bo \in \bO_\a$, and $g\in \cC^{1,2}_{\o,\ab}(\tilde E, E)$, $\tilde g\in \cC^{1,2}_{\tilde \o,\ab}(\tilde E, E)$
define
\bea
\label{dg}
d^{\o,\tilde \o}_{2,\ab}(g, \tilde g) &:=& \|g-\tilde g\|_2 +\|\pa_\o g-\pa_{\tilde\o} \tilde g\|_1\nonumber\\
&&+ \sup_{y\in \tilde E}\Big[d^{\o,\tilde \o}_{\ab}(g(\cd, y), \tilde g(\cd, y)) +  d^{\o,\tilde \o}_{\ab}(\pa_y g(\cd, y), \pa_y \tilde g(\cd, y))\Big].
\eea
\end{defn}

\begin{rem}
\label{rem-space}
{\rm (i) For $g\in \cC^{2,3}_{\o,\ab}(\tilde E, E)$, by \reff{Ito}  we have
\bea
\label{C3rep}
\left.\ba{c}
d g(t,y) = h(t,y) \cd d\bo_t + f(t,y) : d\la\bo\ra_t, \q\mbox{where}\\
h := (\pa_\o g)^* \in \cC^{1,2}_{\bo,\ab,loc}(\tilde E, E^d),  f := D^\o_t g + {1\over 2} \pa_\o h \in \cC^1_{\b, loc}(\tilde E, E^{d\times d}).
\ea\right.
\eea

(ii) In \reff{gC2}, we need only $\|\pa_\o g\|_1$ instead of $\|\pa_\o g\|_{1,\b}$, and in \reff{gC3}, we do not need $\|D^\o_t g\|_{1,\b}$. The latter is particularly convenient because $D^\o_t g$ may not be unique.

(iii) It is clear that $d^\o_{2,\ab}(g,\tilde g) := d^{\o,\o}_{2,\ab}(g, \tilde g) = \|g-\tilde g\|_{2,\bo,\ab}$.
\qed}
\end{rem}

\subsection{Commutativity of $\pa_y$ and path derivatives}


\begin{lem}
\label{lem-commute}
(i) Let $g\in \cC^{2,3}_{\bo,\ab}(\tilde E, E)$. Then $ \pa_{\o y} g = [\pa_{y\o} g]^* \in E^{|\tilde E|\times d}$, namely 
\bea
\label{payg1}
\pa_\o \pa_{y_i} g = \pa_{y_i} \pa_\o g,\q i=1,\cds, |\tilde E|.
\eea 

(ii)  Let $g\in \cC^{3,3}_{\bo,\ab}(\tilde E, E)$.  Then, for appropriate $D^\o_t$ and for each $i=1,\cds, |\tilde E|$, 
\bea
\label{payg2}
\pa^2_{\o\o} \pa_{y_i} g = \pa_{y_i} \pa^2_{\o\o} g &\mbox{and}&  D^\bo_t \pa_{y_i} g =   \pa_{y_i} D^\bo_t g .
\eea 
\end{lem}
\proof  Without loss of generality, we assume $|\tilde E|=1$, namely $\tilde E=\dbR$. Recall \reff{C3rep}.

(i) Fix $y\in \dbR$ and  denote, for $0\neq \D y\in \dbR$,
\beaa
\td \f_t(y) := {\f(t, y+\D y)-\f(t, y)\over \D y}, \q  \f = g, h, f. 
\eeaa   
It is straightforward to check that
\beaa
&\dis \td g_t(y) = \int_0^t \td h_s(y) \cd d\bo_s + \int_0^t \td f_s(y) :d\la\bo\ra_s& \\
 &\dis\td h_t (y)= \int_0^1 \pa_y h(t, y + \l \D y)d\l,\q \td f_t (y)= \int_0^1 \pa_y f(t, y + \l \D y)d\l,&
\eeaa
and thus, as $|\D y|\to 0$,
\beaa
&\dis \bnm\td h(y) - \pa_y h(y)\bnm_{\o,\ab} \le \int_0^1 \bnm\pa_y h(y+\l\D y) - \pa_y h(y)\bnm_{\o,\ab} d\l \to 0,&\\
&\dis \|\td f(y) - \pa_y f(y)\|_{\b} \le \int_0^1 \|\pa_y f(y+\l\D y) - \pa_y f(y)\|_\b d\l \to 0.&
\eeaa
Then it follows from  Lemma \ref{lem-DRI} and Lemma \ref{lem-Young} (ii)   that 
\bea
\label{payg}
\pa_{y} g(t, y) = \int_0^t \pa_{y} h(s, y) \cd d\bo_s+ \int_0^t \pa_{y} f(s, y) : d\la\bo\ra_s.
\eea
This implies \reff{payg1} immediately.

(ii)  Since $h \in \cC^{2,3}_{\bo,\ab}(\tilde E, E^{1\times d})$, by (i) we have $\pa_{y} \pa_\o h = \pa_\o \pa_{y} h$ and thus $\pa_y \pa^2_{\o\o} g = \pa^2_{\o\o}\pa_y g$. Now  applying the convention \reff{Dbot} for $D^\o_t$ on \reff{payg} and by  \reff{C3rep},  we have
\beaa
2D^\bo_t (\pa_y g) &=&(\pa_y f -  {1\over 2}\pa_{\o y}  h) + (\pa_y f -  {1\over 2} \pa_{\o y}  h)^* = \pa_y \Big[ (f -  {1\over 2} \pa_{\o} h) + (f -  {1\over 2} \pa_{\o} h)^*\Big]\\
&=&    (\pa_y f -  {1\over 2} \pa_{y\o} h) + (\pa_y f -  {1\over 2}  \pa_{y\o} h)^* = 2\pa_y D^\bo_t g.
\eeaa
This completes the proof.
\qed


\subsection{Chain rule of path derivatives}

\begin{thm}
\label{thm-chain}
(i) Let $\bo \in \bO_\a$, $\th \in \cC^1_{\o,\a}(\tilde E)$, $g\in \cC^{1,2}_{\o,\ab, loc}(\tilde E, E)$, and $\eta_t := g(t, \th_t)$. Then 
\bea 
\label{chain1}
\eta \in \cC^1_{\o,\a}(E) &\mbox{ with}& \pa_\o \eta_t = (\pa_\o g)(t,\th_t) + \pa_y g(t,\th_t) \otimes \pa_\o \th_t.
\eea

(ii) Assume further that $\th \in \cC^2_{\o,\a}(\tilde E)$ and $g\in  \cC^{2,3}_{\o,\ab, loc}(\tilde E, E)$. Then, for appropriate $D^\bo_t$,   
\bea
\label{chain2}
\eta \in \cC^2_{\o,\a}(E) &\mbox{ with}& D^\bo_t \eta_t = (D^\bo_t g)(t,\th_t) + \pa_y g(t,\th_t) \otimes D^\bo_t \th_t.
\eea
\end{thm}

\begin{rem}
\label{rem-Ventzell} {\rm Similar to \cite{BMZ} Proposition 2.7,  the chain rule of pathwise derivatives is equivalent to the It\^{o}-Ventzell formula, which extends the It\^{o} formula in \cite{FH} Proposition 5.6. Indeed, note that $\th\in \cC^2_{\o,\a}(\tilde E)$ takes the form:
\bea
\label{threp}
d\th_t = a_t \cd d\bo_t + b_t : d\la\bo\ra_t &\mbox{where}&a := (\pa_\o \th)^*, ~ b := D^\bo_t\th + {1\over 2} \pa_\o a.
\eea
Recall \reff{C3rep} again.  It follows from Lemma \ref{lem-commute} (i) that $\pa_\o \pa_y g = (\pa_y h)^*$. Then, noticing that $h \in \cC^{1,2}_{\o,\ab, loc}(\tilde E, E^d)$, $\pa_y g \in \cC^{1,2}_{\o,\ab, loc}(\tilde E, E^{1\times |\tilde E|})$, by applying \reff{chain1} several times and by  \reff{chain2}, we have
\beaa
\pa_\o \eta_t &=& h^*(t,\th_t) + \pa_y g(t,\th_t)\otimes a^*_t,\\
\pa^2_{\o\o} \eta_t &=& \pa_\o[h(t,\th_t) + \pa_y g(t,\th_t)\otimes a_t] \\
& = &\Big[\pa_\o h + \pa_y h \otimes a^* + (\pa_y h \otimes a^*)^* + \pa^2_{yy} g \otimes_2 [a, a] + \pa_y g \otimes \pa_\o a\Big](t, \th_t);\\
 D^\bo_t \eta_t &=& {1\over 2}\Big[ [(f - {1\over 2}\pa_\o h)+(f-{1\over 2}\pa_\o h)^*] + \pa_y g \otimes [(b-{1\over 2}\pa_\o a) + (b-{1\over 2}\pa_\o a)^*]\Big](t,\th_t).
\eeaa
 This, together with \reff{Ito} and the symmetry of $\la\bo\ra$, implies: 
\bea
\label{Ventzell}
d[g(t, \th_t)] 
&=& \Big[ h(t,\th_t) +  \pa_y g(t,\th_t) \otimes  a_t\Big]\cd d\bo_t \\
&&+ \Big[ f +  \pa_y g \otimes  b_t + {1\over 2}  \pa^2_{yy} g\otimes_2  [a_t, a_t] +    \pa_{y} h\otimes  a^*_t \Big](t, \th_t) :  d\la \bo\ra_t,\nonumber
 \eea
which we call the  pathwise  It\^{o}-Ventzell formula. 
\qed}
\end{rem}

\no{\bf Proof of Theorem \ref{thm-chain}.} (i) For $(s, t)\in \dbT^2$, we have
\bea
\label{Deta}
\eta_{s,t} &=&  g(t, \th_t) - g(s, \th_s) = g(t, \th_t) - g(s, \th_t) + g(s, \th_t) - g(s, \th_s) \\
&=&  [\pa_\o g](s, \th_t) \o_{s,t} + R^{\o, g(\cd, \th_t)}_{s,t} + \int_0^1 \pa_y g(s, \th_s + \l \th_{s,t}) d\l\otimes \th_{s,t} \nonumber\\
&=& \Big[(\pa_\o g)(s,\th_s) + \pa_y g(s,\th_s) \otimes \pa_\o \th_s\Big] \o_{s,t} + R^{\o, \eta}_{s,t},\nonumber
\eea
where
\beaa
&&R^{\o, \eta}_{s,t} :=\Big[[\pa_\o g](s, \th_t)-[\pa_\o g](s,\th_s)\Big] \o_{s,t} + R^{\o, g(\cd, \th_t)}_{s,t} \\
&&\q+ \int_0^1 [\pa_y g(s, \th_s + \l \th_{s,t}) -\pa_y g(s,\th_s)] d\l\otimes \pa_\o \th_s \o_{s,t} +   \int_0^1 \pa_y g(s, \th_s + \l \th_{s,t}) d\l\otimes R^{\o, \th}_{s,t}.
\eeaa
Then clearly
\bea
\label{Reta}
\|R^{\o,\eta}\|_{\a+\b} \le \|g\|_{2,\o,\ab}\Big[\|\th\|_\b \|\o\|_\a + 1+ \|\th\|_\b \|\pa_\o\th\|_\infty \|\o\|_\a + \|\th\|_{\o, \ab}\Big] < \infty.
\eea
Moreover, under our conditions it is clear that $(\pa_\o g)(t,\th_t) + \pa_y g(t,\th_t) \otimes \pa_\o \th_t$ is H\"{o}lder-$\b$-continuous. This proves \reff{chain1}.

(ii) Recall \reff{C3rep} and \reff{threp}. By reversing the arguments in Remark \ref{rem-Ventzell}, it suffices to prove \reff{Ventzell}. Denote $\d:= t-s$. Recall the first line of \reff{Deta} and note that 
\beaa
&\th_{s, t} =  a_s \cd \o_{s,t} + \pa_\o a_s: \uo_{s,t} + b_s : \la\bo\ra_{s,t} + O(\d^{2\a+\b});&\\
&g(t,y)-g(s,y) = h(s,y) \cd \o_{s,t} + \pa_\o h(s,y): \uo_{s,t} + f(s,y) : \la\bo\ra_{s,t}  + O(\d^{2\a+\b})&
\eeaa
Then, by the standard Taylor expansion and applying Lemma \ref{lem-commute} (i) on $g$, we have
\beaa
&&g(t,\th_t) - g(t,\th_s) = \pa_y g(t, \th_s) \otimes  \th_{s, t} + {1\over 2} \pa^2_{yy} g(t, \th_s) \otimes_2  [\th_{s,t}, \th_{s,t}] + O(\d^{3\a}) \\
&&\qq= \Big[\pa_y g(s, \th_s) + \pa_y h(s, \th_s) \cd \o_{s,t}\Big]  \otimes  \th_{s, t} +  {1\over 2} \pa^2_{yy} g(s, \th_s) \otimes_2 [ \th_{s,t}, \th_{s,t}] + O(\d^{2\a+\b});\\
&&g(t,\th_s) - g(s,\th_s) =h(s,\th_s) \cd \o_{s,t} + [\pa_\o h](s,\th_s): \uo_{s,t} + f(s,\th_s) : \la\bo\ra_{s,t} + O(\d^{2\a+\b}).
\eeaa
On the other hand,
\beaa
&&\int_s^t [ h(r, \th_r)  + \pa_y g(r, \th_r) \otimes   a_r] \cd d\bo_r\\
&&\q = [ h(s, \th_s)  + \pa_y g(s, \th_s) \otimes   a_s]\cd \o_{s,t} + \pa_\o [ h(s, \th_s)  + \pa_y g(s, \th_s) \otimes   a_s] : \uo_{s,t} + O(\d^{2\a+\b});\\
&&\int_s^t [f(r,\th_r) + \pa_y g(r,\th_r)\otimes b_r] : d\la\bo\ra_r = [f(s,\th_s) + \pa_y g(s,\th_s)\otimes b_s] : \la\bo\ra_{s,t} +  O(\d^{2\a+\b}).
\eeaa
By Lemma \ref{lem-commute} (i) we have $\pa_{\o y} g = [\pa_{y \o} g]^*=\pa_y h^*$. Then it follows from \reff{chain1} that
\bea
\label{paoh}
&&\pa_\o [ h(s, \th_s)  + \pa_y g(s, \th_s) \otimes   a_s]\\
&=& \Big[\pa_\o h  + \pa_y h \otimes  a^*_s + \pa_y h^* \otimes   a_s+ \pa^2_{yy} g \otimes_2 [a_s^*, a^*_s] + \pa_y g \otimes  \pa_\o a_s ](s, \th_s).\nonumber
\eea
Noting that $\o_{s,t} = O(\d^\a)$, $\uo_{s,t} = O(\d^{2\a})$, and $\la\bo\ra_{s,t} = O(\d^{2\a})$, then we have
\beaa
&& \eta_{s,t}  - \int_s^t   [ h(r, \th_r)  + \pa_y g(r, \th_r) \otimes   a_r] \cd d\bo_r -  \int_s^t [f(r,\th_r) + \pa_y g(r,\th_r)\otimes b_r] : d\la\bo\ra_r\\
&=&  \Big[ [ \pa_y h(s, \th_s) \cd \o_{s,t}]  \otimes  [a_s \cd \o_{s,t}]+  {1\over 2} \pa^2_{yy} g(t, \th_s) \otimes_2 [(a_s \cd \o_{s,t})^*, (a_s \cd \o_{s,t})^*] \\
&&-  \Big[  \pa_y h(s,\th_s) \otimes  a^*_s + [\pa_y h(s,\th_s) \otimes  a^*_s]^*  + \pa^2_{yy} g(s,\th_s) \otimes_2  [a^*_s, a^*_s]]\Big] : \uo_{s,t} +  O(\d^{2\a+\b})\\
&=& \Big[ {1\over 2} \pa^2_{yy} g(t, \th_s) \otimes_2 [\pa_\o \th_s, \pa_\o \th_s] + \pa_y h (s,\th_s) \otimes \pa_\o \th_s\Big] : \la\bo\ra_{s,t} +  O(\d^{2\a+\b})
\eeaa
This proves \reff{Ventzell}, and hence \reff{chain2}.
\qed

\subsection{Some estimates} 
In this subsection we provide some estimates for $\eta = g(t,\th_t)$, which will be crucial for studying rough differential equations in next section. These results correspond to  \cite{FH} Lemma 7.3 and Theorem 7.5, where $g$ does not depend on $t$. 
\begin{lem}
\label{lem-g}
(i) Let $\bo\in\bO_\a$,  $\th\in \cC^1_{\o,\ab}(E)$, $g \in \cC^{1,2}_{\o,\ab}(\tilde E,E)$, $\eta_t := g(t, \th_t)$, and denote 
\beaa
M_1 := \|\bo\|_\a+\bnm \th\bnm_{\o,\ab}.
\eeaa
Then for any $T_0>0$ and any $T\le T_0$,  there exists a constant $C_{\ab, M_1, T_0}$, depending only on $\ab$, $M_1$, $T_0$, and $|E|$, $|\tilde E|$, such that
\bea
\label{etaest}
\|\eta\|_{\o,\ab} &\le& C_{\ab, M_1, T_0} \|g\|_{2,\o, \ab}.
\eea

(ii) Assume further that $g \in \cC^{2,3}_{\o,\ab}(\tilde E, E)$, and  $(\tilde \bo, \tilde \th, \tilde g, \tilde \eta)$ satisfy the same conditions. Denote $\D \f := \tilde \f-\f$ for appropriate $\f$, and 
\beaa
 M_2:=  \bnm \th\bnm_{\o,\ab} +  \bnm\tilde \th\bnm_{\tilde\o,\ab} +\|\bo\|_\a + \|\tilde\bo\|_\a +  \|g\|_{3,\o,\ab} +  \|\tilde g\|_{3,\tilde \o,\ab}.
 \eeaa
Then, for any $T\le T_0$ as in (i), there exists a constant $ C_{\ab, M_2, T_0}$ such that
\bea
\label{Detaest}
d^{\o,\tilde\o}_\ab(\eta, \tilde \eta) \le C_{\ab, M_2, T_0} \Big[d^{\o, \tilde\o}_{2,\ab}(g, \tilde g) + \bd^{\o,\tilde\o}_\ab(\th, \tilde\th) + |\D\th_0|+\|\D\bo\|_\a\Big].
\eea
\end{lem}
\proof (i) First, by \reff{oinfty} and \reff{tha} we have $\|\pa_\o \th\|_\infty+\|\th\|_\b \le C$. By the first line of \reff{Deta} it is clear that
\bea
\label{gest}
\|\eta\|_\b &\le& C\Big[\|g\|_{0,\b} + \|g\|_{1}\Big].
\eea
Next, recall \reff{chain1} and note that 
\beaa
|\pa_\o \eta_{s,t}| \le |\pa_\o g(t, \th_t) - \pa_\o g(s, \th_s)| + |\pa_y g(t, \th_t) - \pa_y g(s, \th_s)| |\pa_\o \th_t| + |\pa_y g(s, \th_s)||\pa_\o \th_{s,t}|.
\eeaa
Applying \reff{gest} on $\pa_\o g$ and $\pa_y g$ we obtain $\|\pa_\o \eta\|_\b\le C\|g\|_{2,\o, \ab}$. Moreover, 
by \reff{Reta} we have $\|R^{\o,\eta}\|_{\a+\b} \le C\|g\|_{2,\o, \ab}$. Putting together we prove  \reff{etaest}.

(ii) First, note that
\beaa
\D\eta_{s,t} &=& \tilde g(t, \tilde \th_t) - g(t, \th_t) - \tilde g(s, \tilde \th_s) + g(s, \th_s)\\
&=& [\D g (t, \tilde \th_t) - \D g (s, \tilde \th_s) ]  + \int_0^1 \pa_y g(s, \th_s+ \l \D \th_s) d\l \otimes  \D \th_{s,t}  \\
&&+  \int_0^1 [\pa_y g(t, \th_t + \l \D\th_t) -\pa_y g(s,  \th_s + \l\D\th_s) ] d\l  \otimes \D \th_t.
\eeaa
Apply \reff{gest} on $\D g$ and $\pa_y g$, we obtain
\beaa
\|\D \eta\|_\b &\le&C\Big[ \|\D g\|_{0,\b}+ \|\D g\|_{1}  + \|\D \th\|_\b +  |\D\th_0|\Big]
\eeaa
Note that $\th_{s,t} = \pa_\o \th_s \o_{s,t} + R^{\o,\th}_{s,t}$, and similarly for $\tilde \th$. Then, by \reff{oinfty}, 
\bea
\label{tha2}
\|\D \th\|_\b &\le& \| \pa_{\tilde \o}\tilde \th-\pa_\o \th\|_\infty \|\tilde \o\|_\b + \|\pa_\o \th\|_\infty \|\D \o\|_\b + \|R^{\tilde\o,\tilde\th}-R^{\o,\th}\|_\b\nonumber\\
& \le& C\Big[\bd^{\o,\tilde \o}_\a(\th, \tilde\th)] + \|\D \o\|_\a\Big].
\eea
Thus
\bea
\label{Detaest1}
\|\D \eta\|_\b&\le& C\Big[\|\D g\|_{0,\b}+ \|\D g\|_{1}  + |\D \th_0|+ \bd^{\o,\tilde \o}_\a(\th, \tilde\th)+\|\D \o\|_\a\Big].
\eea
We shall emphasize that the above $C$ depends on $ \|g\|_{2,\o,\ab} +  \|\tilde g\|_{2,\tilde \o,\ab}$, not  $\|g\|_{3,\o,\ab} +  \|\tilde g\|_{3,\tilde \o,\ab}$.

Next, note that
\beaa
&&\pa_{\tilde \o}\tilde \eta_t -\pa_\o \eta_t = [\pa_{\tilde \o} \tilde g(t, \tilde \th_t) -\pa_\o g(t, \th_t)] + [\pa_y \tilde g (t, \tilde \th_t) - \pa_y g (t,  \th_t)] \otimes \pa_{\tilde \o} \tilde \th_t\\
&& \qq\qq\qq + \pa_y g(t, \th_t)\otimes [\pa_{\tilde \o}\tilde  \th_t -\pa_\o\th_t]. \\
&&[\pa_{\tilde \o}\tilde \eta -\pa_\o \eta]_{s,t} =  [\pa_{\tilde\o} \tilde g(\cd, \tilde \th_\cd) - \pa_\o g(\cd, \th_\cd)]_{s,t} +  [\pa_y \tilde g (\cd, \tilde \th_\cd) - \pa_y g (\cd,  \th_\cd)]_{s,t} \otimes \pa_\o \tilde \th_t \\
&&\qq\qq\qq\q+  [\pa_y \D g (s, \tilde \th_s) + \pa_y g(s, \tilde \th_s) - \pa_y g (s,  \th_s)] \otimes \pa_\o \tilde \th_{s,t} \\
&&\qq\qq\qq\q + [\pa_y g (\cd, \th_\cd) ]_{s,t} \otimes \D \pa_\o \th_t + \pa_yg(s, \th_s) \otimes \D \pa_\o \th_{s,t}.
\eeaa 
Apply \reff{Detaest1} on $\pa_\o g$ and $\pa_y g$, and \reff{gest} on $\pa_y g$, we obtain from \reff{tha2} that
\bea
\label{Deta'est}
\|\D \pa_\o \eta\|_\a &\le& C\Big[d^{\o, \tilde\o}_{2,\a}(g, \tilde g) +|\D \th_0|+\bd^{\o,\tilde \o}_\a(\th, \tilde\th)+\|\D \o\|_\a\Big]
\eea

Finally, recall \reff{Reta} and note that 
\beaa
&& R^{\tilde\o, \tilde g(\cd, \tilde y)}_{s, t} - R^{\o, g(\cd, y)}_{s, t} \\
&=&R^{\tilde\o, \tilde g(\cd, \tilde y)}_{s, t} - R^{\o, g(\cd, \tilde y)}_{s, t} +  \Big[[g(\cd, \tilde y)\_{s,t}  - \pa_\o g(s, \tilde y)  \o_{s, t}\Big] - \Big[ [g(\cd, y)]_{s,t}  - \pa_\o g(s, y)  \o_{s, t}\Big]\\
&=&R^{\tilde\o, \tilde g(\cd, \tilde y)}_{s, t} - R^{\o, g(\cd, \tilde y)}_{s, t} + \int_0^1 R^{\o, \pa_y g(\cd, y +\l \D y)}_{s, t}d\l \otimes \D y,
\eeaa
one can obtain the desired estimate for $\|R^{\tilde\o, \tilde\eta}-R^{\o, \eta}\|_{\a+\b}$ straightforwardly. This, together with  \reff{Deta'est}, completes the proof.  
\qed

Moreover, we have the following simpler results whose proof is omitted. 
\begin{lem}
\label{lem-Youngf} (i) Let $\th \in \O_\b(E)$, $f\in \cC^1_\b(\tilde E, E)$, and $\eta_t := f(t, \th_t)$.  Then $\eta \in \O_\b(E)$ and 
\bea
\label{feta}
\|\eta\|_\b \le \|f\|_{0,\b} + \|f\|_{1} \|\th\|_\b \le \|f\|_{1,\b}[1+\|\th\|_\b].
\eea
(ii) Let $\th,\tilde\th \in \O_\b(E)$, $f, \tilde f\in \cC^2_\b(\tilde E,  E)$, and $\eta_t := f(t, \th_t)$, $\tilde \eta:= \tilde f(t, \tilde \th_t)$.  Then
\bea
\label{Dfeta}
\|\tilde \eta - \eta \|_{\b} \le [1+\|\th\|_\b+ \|\tilde\th\|_\b] \Big[\|\tilde f - f\|_{1,\b} + \|f\|_{2} [|\tilde\th_0-\th_0|+\|\tilde\th - \th\|_\b]\Big]. 
\eea
\end{lem}

\section{Rough Differential Equations}
\label{sect-RDE}
\setcounter{equation}{0}
In this section we study rough path differential equations  with coefficients less regular in the time variable $t$, motivated from our study of stochastic differential equations  with random coefficients in next section. 
Let $\bo\in \bO_\a$, $g\in  \cC^{2,3}_{\o,\ab}(E, E^d)$, $f\in \cC^2_\b(E, E^{d\times d})$, and $y_0 \in E$. Consider the following RDE:
\bea
\label{RDE}
\th_t = y_0 +\int_0^t g(s, \th_s) \cd d\bo_s + \int_0^t f(s, \th_s): d\la \bo\ra_s , \q t\in \dbT.
\eea
Our goal is to find solution $\th \in \cC^1_{\o,\ab}(E)$. By Theorem \ref{thm-chain} and Lemma \ref{lem-Youngf}, in this case $g(\cd, \th)\in \cC^1_{\o,\ab}(E^d)$, $f(\cd, \th)\in \O_\b(E^{d\times d})$, and thus the right side of \reff{RDE} is well defined.

\begin{rem}
\label{rem-sol}
{\rm When $\th\in  \cC^1_{\o,\ab}(E)$ is a solution, clearly $\pa_\o \th_t = g(t,\th_t)$, then by Theorem \ref{thm-chain} (i) it is clear that $\th \in \cC^2_{\o, \ab}(E)$. So a solution to RDE \reff{RDE} is automatically in  $\cC^2_{\o, \ab}(E)$. We shall use this fact without mentioning it.
\qed}
\end{rem}

In standard rough path theory the generator $g$ of  RDE \reff{RDE} is independent of $t$. In Lejay and Victoir \cite{LV}, $g$ may depend on $t$, but is required to be H\"{o}lder-$(1-\a)$ continuous, which is violated for $g\in  \cC^{2,3}_{\o,\ab}(E, E^d)$ (since $\a<{1\over 2}$). This relaxation of regularity in $t$ is crucial for studying SDEs and SPDEs with random coefficients, see Remark \ref{rem-SDE} below. We also refer to Gubinelli, Tindel and Torrecilla \cite{GTT} for some discussion along this direction.

\begin{thm}
\label{thm-RDEwell}
Let $\bo\in \bO_\a$, $g\in \cC^{2,3}_{\o,\ab}(E, E^d)$, $f\in \cC^2_\b(E, E^{d\times d})$, and $y_0 \in E$. Then RDE \reff{RDE} has a unique solution  $\th \in \cC^2_{\o,\ab}(E)$. Moreover, there exists a constant $C_\ab$, depending only on $\ab$, $d$, $|E|$, $T$, $\|f\|_{2,\b}$, $\|g\|_{3, \o,  \ab}$, and $\|\bo\|_\a$, such that 
\bea
\label{RDEbound}
\|\th\|_\a+ \| \th\|_{\o, \ab} \le C_{\ab}.
\eea
\end{thm}
\proof We proceed in three steps.

{\it Step 1.} Denote $M:= [\|\pa_\o g\|_{0}+\|g\|_{1}^2]  \|\bo\|_\a+\|f\|_{0} \|\bo\|_\a[2+\|\bo\|_\a]$ and 
\bea
\label{cAb}
\cA_\ab := \Big\{\th\in \cC^1_{\o,\ab}(E): \th_0 = y_0, \pa_\o \th_0 = g^*(0, y_0), \|\th\|_{\o, \ab} \le M+1\Big\}, 
\eea
equipped with the norm $\|\cd\|_{\o, \ab}$. Note that $\cA_\ab$ contains $\th_t := y_0 + g(0, y_0)\cd \o_{0,t}$ and thus  is not empty.  Define a mapping $\Phi$ on $\cA_\ab$:
\beaa
\Phi(\th) :=\Th ~~\mbox{where}~ \Th_t :=  y_0 +\Th^1_t + \Th^2_t := y_0+ \int_0^t g(s,\th_s)\cd d\bo_s +\int_0^t f(s, \th_s):d\la\bo\ra_s.
\eeaa
We show that, there exists $0<\d\le 1$, which depends on $\ab$, $d$, $|E|$, $T$, $\|f\|_{2,\b}$, $\|g\|_{3,\o, \ab}$,  and $\|\bo\|_\a$, but not on $y_0$, such that whenever $T\le \d$, $\Phi$ is a contraction mapping on $\cA_\ab$. One can easily check that $\cA_\ab$ is complete under $d^{\o,\o}_\ab$, then $\Phi$ has a unique fixed point $\th \in \cA_\ab$ which is clearly the unique solution of RDE \reff{RDE}.

To prove that $\Phi$ is a contraction mapping, let $C$ denote a generic constant which depends only on the above parameters, but not on $y_0$.
We first show that $\Phi(\th)\in \cA_\ab$ for all $\th \in \cA_\ab$. Indeed, clearly $\Th_0=y_0$ and $\pa_\o \th_0= g^*(0, y_0)$. For any $\th\in \cA_\ab$, denote $\eta_t := g(t, \th_t)$. Applying Lemma  \ref{lem-g}  and then Lemma \ref{lem-RI}, we have, 
\beaa
&\|\eta\|_{\o,\ab} \le C,\q \|\pa_\o \eta_0| \le  \|\pa_\o g\|_{0}+\|\pa_\o g\|_{1}^2,\q \mbox{and thus}&\\
& \|\Th^1\|_{\o, \ab} \le \|\bo\|_\a|\pa_\o \eta_0| + C\d^\a[1+\|\bo\|_\a]\|\eta\|_{\o,\ab} \le   [\|\pa_\o g\|_{0}+\|g\|_{1}^2] \|\bo\|_\a+ C \d^\a.&
\eeaa
Similarly, It follows from Lemmas \ref{lem-Young}  and   \ref{lem-Youngf} (i)  that
\beaa
 \|\Th^2\|_{\o,\ab} = \|\Th^2\|_{\a+\b} \le \|f\|_{0} \|\bo\|_\a[2+\|\bo\|_\a] + C\d^\a,\\
   \mbox{and thus}~
\|\Th\|_{\o,\ab} \le \|\Th^1\|_{\o,\ab} +  \|\Th^2\|_{\o, \ab} \le M+ C \d^\a.
\eeaa
Set $\d$ small enough we have $\|\Th\|_{\o,\ab}  \le M+1$. That is, $\Th\in \cA_\ab$.

Next, let $\tilde \th \in \cA_\ab$ and denote $\tilde\Th, \tilde\Th^1, \tilde\Th^2, \tilde\eta$ in obvious sense. Let $\D \f := \tilde \f - \f$ for appropriate $\f$. Recall \reff{tha2} we see that
\bea
\label{Dthinfty}
\|\D \th\|_\infty \le C\d^\b \|\D \th\|_\b \le C\d^\b \|\D\th\|_\ab.
\eea
Then, applying Lemmas  \ref{lem-DRI},   \ref{lem-g} (ii),  \ref{lem-Young} (ii), and  \ref{lem-Youngf} (ii), we have
\beaa
&\|\D\Th^1\|_{\o,\ab} \le C\d^\a  \|\D\eta\|_{\o,\ab} \le C\d^\a \|\D\th\|_{\o,\ab},\q \|\D\Th^2\|_{\a+\b} \le  C\d^{\a} \|\D\th\|_\b,&\\
&\mbox{and thus}\q 
\|\D\Th\|_{\o,\ab} \le C\d^\a \|\D\th\|_{\o,\ab}.&
\eeaa
Set $\d$ be small enough such that $C\d^{\a} \le {1\over 2}$, then $\Phi$ is a contraction mapping.

{\it Step 2.}   We now prove the result for general $T$. Let  $\d$ be the constant in Step 1. Let $0=t_0<\cds<t_n = T$ such that $t_{i+1}-t_i\le \d$, $i=0,\cds, n-1$. We may solve the RDE over each interval $[t_i, t_{i+1}]$ with initial condition $(\th_{t_i}, g(t_i, \th_{t_i}))$, which is obtained  from the previous step by considering the RDE  on $[t_{i-1}, t_i]$, and thus we obtain the unique solution over the whole interval $[0, T]$. 

{\it Step 3.} We now estimate $\|\th\|_{\o,\ab}$. First,  when $T \le \d$ for the constant $\d = \d_{\ab}$ in {\it Step 1}, we have $\th\in \cA_\ab$ and thus $\|\th\|_\b \le M+1$. In particular, this implies that
\beaa
|\pa_\o \th_{s,t}| \le (M+1)(t-s)^\b,\q |R^{\o, \th}_{s,t}|\le (M+1)(t-s)^{\a+\b},\q\mbox{whenever}~t-s\le \d.
\eeaa
Now for arbitrary $s, t$, let  $k := [{t-s\over \d}] + 1$ be the smallest integer greater than ${t-s\over \d}$,  and $t_i := s + {i\over k}(t-s)$, $i=0,\cds, k$. Then
\beaa
|\pa_\o \th_{s,t}| \le \sum_{i=0}^{k-1} |\pa_\o \th_{t_i, t_{i+1}}| \le k ({t-s\over k})^\b = k^{1-\b} (t-s)^\b \le (\d^{-1}+1)^{1-\b} (t-s)^\b.
\eeaa 
Thus $\|\pa_\o \th\|_\b \le (\d^{-1}+1)^{1-\b}$. Similarly we may prove $\|R^{\o, \th}\|_{\a+\b} \le  (\d^{-1}+1)^{1-\a-\b}$.

Finally, note that $\|\pa_\o \th\|_\infty \le C$, it is clear that $\|\th\|_\a \le \|\pa_\o \th\|_\infty \|\o\|_\a + \|R^{\o,\th}\|_\a \le C$. 
\qed

\bs

We next study the stability of RDEs.
\begin{thm}
\label{thm-RDEstability}
Let $(y_0, \bo, f, g)$ and $(\tilde y_0, \tilde\bo, \tilde f, \tilde g)$ be as in Theorem \ref{thm-RDEwell}, and $\th$, $\tilde\th$ be the corresponding solution of the RDE. 
Then there exists a constant $C_\ab$, depending only on  $\ab$, $d$, $|E|$, $T$, $\|f\|_{2,\b}$, $\|\tilde f\|_{2,\b}$, $\|g\|_{3, \o, \ab}$, $\|\tilde g\|_{3, \tilde\o, \ab}$, and $\|\bo\|_\a$, $\|\tilde \bo\|_{\a}$, such that, denoting $\D \f := \f-\tilde \f$ for appropriate $\f$,
\bea
\label{RDEStability}
d^{\o,\tilde\o}_\ab(\th, \tilde\th) \le C_{\ab}[\D I_\ab + |\D y_0|]~\mbox{where}~ \D I_\ab:=d^{\o, \tilde \o}_{2,\ab}(g, \tilde g) + \|\D f\|_{1,\b} + \|\D\bo\|_\a.
\eea
\end{thm}
\proof First assume $T \le \d$ for some constant $\d>0$ small enough. Use the notations in Step 1 of Theorem  \ref{thm-RDEwell}. Applying Lemma \ref{lem-g} (i) and \reff{RDEbound} we see that $|\pa_{\tilde\o}\tilde \eta_0|+ \|\tilde \eta\|_{\o,\b} \le C$. Then, it follows from
 Lemmas \ref{lem-DRI}  and \ref{lem-g} (ii) that
\beaa
d^{\o,\tilde\o}_\ab(\Th^1, \tilde \Th^1) \le C\Big[\d^\a  d^{\o,\tilde\o}_\ab(\eta, \tilde \eta) + d_\a(\bo, \tilde \bo) + |\eta'_0 - \tilde \eta'_0|\Big]\le C\Big[\d^\a d^{\o,\tilde\o}_\ab(\th, \tilde \th) + \D I_\ab+|\D y_0|\Big].
\eeaa
Similarly, by Lemmas \ref{lem-Young} and  \ref{lem-Youngf},    we have
\beaa
\|\D\Th^2\|_{\a+\b} \le C\Big[\d^\a \|\D\th\|_\b  + \D I_\ab+|\D y_0|\Big].
\eeaa
Putting together we get
\beaa
d^{\o,\tilde\o}_\ab(\th, \tilde \th) = d^{\o,\tilde\o}_\ab(\Th, \tilde \Th) \le C\Big[\d^\a d^{\o,\tilde\o}_\ab(\th, \tilde \th) + \D I_\ab+|\D y_0|\Big].
\eeaa
Set $\d$ be small enough such that $C\d^{\a} \le {1\over 2}$, we obtain $d^{\o,\tilde\o}_\ab(\th, \tilde \th) \le C[\D I_\ab+|\D y_0|]$.

Now for general $T$, let $k:= [{T\over \d}]+1$ be the smallest integer greater than ${T\over \d}$ and  $t_i := {i\over k} T$, $i=0,\cds, k$. Denote
\beaa
\D J_i := \sup_{t_i \le s < t\le t_{i+1}} \Big[{|\D\pa_\o \th_{s,t}|\over (t-s)^\b} + {|R^{\tilde\o, \tilde\th}_{s,t}-R^{\o, \th}_{s,t}|\over (t-s)^{\a+\b}}\Big],\q i=0,\cds, k-1.
\eeaa
By the above arguments we have $\D J_i \le C[\D I_\ab + |\D \th_{t_i}|]$. Then, applying \reff{tha2} on $[t_i, t_{i+1}]$ and noting that $\pa_\o \th_{t_i}=g(t_i, \th_{t_i})$ and $\pa_\o \tilde \th_{t_i}=\tilde g(t_i, \tilde \th_{t_i})$ are bounded, we have
\beaa
 |\D \th_{t_{i+1}}| \le |\D \th_{t_i}| + |\D \th_{t_i, t_{i+1}}| \le  |\D \th_{t_i}| + \D J_i + C[|\D \pa_\o \th_{t_i}| + \|\D \o\|_\a] \le   C[\D I_\ab + |\D \th_{t_i}|].
\eeaa
By induction we get 
\beaa
\max_{0\le i\le k} |\D \th_{t_i}| \le C[\D I_\ab + |\D y_0|], &\mbox{and thus}&\max_{0\le i\le k}\D J_i \le C[\D I_\ab + |\D y_0|].
\eeaa
 Now following the arguments  in Theorem \ref{thm-RDEwell} Step 3 we can prove the desired estimate.
\qed

\begin{rem}
\label{rem-RDE}
{\rm 

(i) The uniqueness of RDE solutions do not depend on boundedness of $g$, $\pa_\o g$, and $f$. Indeed, let $\th$ and $\tilde\th$ be two solutions. Notice that any element of $\cC^1_{\o,\a}(E)$ is bounded, and thus we may denote $M_0 := \|\th\|_\infty + \|\tilde\th\|_\infty<\infty$.  One can see that all the arguments in Theorem \ref{thm-RDEwell} remain valid if we replace the $\sup_{y\in E}$ in \reff{gk} with $\sup_{y\in E, |y|\le M_0}$, while the latter is always bounded for $g$, $\pa_\o g$, and $f$. 

(ii) If we do not assume boundedness of $g$, $\pa_\o g$, and $f$, in general we can only obtain the local existence, namely the solution exists when $T$ is small. However, if we can construct a solution for large $T$, as we will see for linear RDEs, then by (ii) above this solution is the unique solution.  
\qed}
\end{rem}

\subsection{Linear RDE}
\label{sect-LRDE}

Now consider RDE \reff{RDE} with
\bea
\label{linear} 
\left.\ba{c}
g(t, y) = a_t \otimes y +b_t,~ f(t,y) = \l_t \otimes y + l_t,~  ,~\mbox{where}~\\
 y\in E, a\in \cC^2_{\o,\ab}(E^{d\times |E|}), b \in \cC^1_{\o,\ab}(E^d), \l \in \O_\b(E^{d\times d\times |E|}), l \in \O_\b(E^{d\times d}).
\ea\right.
\eea
We remark that the above $f$ and $g$ are not bounded and thus we cannot apply Theorem \ref{thm-RDEwell} directly. In Friz and Victoir \cite{FV}, some a priori estimate is provided for linear RDEs and then the global existence follows from the arguments of Theorem \ref{thm-RDEwell}, by replacing the $\sup_{y\in E}$ in \reff{gk}  with the supremum over the a priori bound of the solution, as illustrated in Remark \ref{rem-RDE} (ii). At below, we shall construct a solution semi-explicitly.  When $|E|=1$, we have an explicit representation in the spirit of Feyman-Kac formula in stochastic analysis literature, see \reff{linearrep} below. However, the formula fails in multidimensional case due to the noncommutativity of matrices. Our main idea is to introduce a decoupling strategy, by using the local solution of certain Riccati type of RDEs, so as to reduce the dimension of $E$.  To our best knowledge, such a construction is new even for  multidimensional  linear SDEs. 

\begin{thm}
\label{thm-RDElinear} The linear RDE \reff{RDE} with \reff{linear} has a unique solution.
\end{thm}
\proof  If $b \in \cC^2_{\o,\ab}(E^d)$, under \reff{linear}  it is straightforward to check that  $g\in \cC^{2,3}_{\o,\ab,loc}(E, E^d)$ and $f\in \cC^2_{\b,loc}(E, E^{d\times d})$, and thus the uniqueness follows from Theorem \ref{thm-RDEwell} and  Remark \ref{rem-RDE} (ii). 
However, in the linear case, by going through the arguments of Theorem \ref{thm-RDEwell} we can easily see that it is enough to assume the  weaker condition $b \in \cC^1_{\o,\ab}(E^d)$. 
We shall construct the solution and thus obtain the existence via induction on $|E|$. 

{\it Step 1.} We first assume  $|E|=1$, namely $E=\dbR$.  Applying Theorem \ref{thm-chain} and Remark \ref{rem-Ventzell} we may verify directly that the following provides a representation of the solution:
\bea
\label{linearrep}
&\dis\th_t = \G^{-1}_t \Big[\th_0 + \int_0^t \G_s b_s \cd d\bo_s + \int_0^t \G_s \big[l_s - a_s b^*_s\big] : d\la \bo\ra_s\Big],&\\
&\dis \mbox{where} ~\G_t:= \exp\Big( -\int_0^t a_s \cd d\bo_s + \int_0^t\big[ {1\over 2} a_s a^*_s-l_s\big] : d\la \bo\ra_s\Big). &\nonumber
\eea

{\it Step 2.} In order to show the induction idea clearly, we present the case $|E|=2$ in details. With the notations in obvious sense, the linear RDE becomes
\bea
\label{linear2}
\left.\ba{c}
d \th^1_t = [a^{11}_t \th^1_t + a^{12}_t \th^2_t + b^1_t] \cd d\bo_t + [\l^{11}_t \th^1_t + \l^{12}_t \th^2_t + l^1_t] : d\la\bo\ra_t;\\
d \th^2_t = [a^{21}_t \th^1_t + a^{22}_t \th^2_t + b^2_t] \cd d\bo_t + [\l^{21}_t \th^1_t + \l^{22}_t \th^2_t + l^2_t] : d\la\bo\ra_t.
\ea\right.
\eea
Clearly, if the system is decoupled, for example if $a^{12} = 0$ and $\l^{12}=0$, one can easily solve the system by first solving for $\th^1$ and then solving for $\th^2$. In the general case, we  introduce a decoupling strategy as follows. Consider an auxiliary RDE:
\bea
\label{barG}
d \ol \G_t = \ol a_t \cd d\bo_t + \ol \l_t : d\la\bo\ra_t.
\eea
where $\ol a, \ol \l$ will be specified later. Denote $\ol \th_t := \th^2_t + \ol \G_t \th^1_t$. Then, applying the It\^{o}-Ventzell formula \reff{Ventzell} we have
\bea
\label{barth}
&&d\ol \th_t = \Big[[a^{22}_t \th^2_t + a^{21}_t \th^1_t + b^2_t] + \ol\G_t [a^{12}_t \th^2_t + a^{11}_t \th^1_t + b^1_t] + \ol a_t \th^1_t \Big] \cd d\bo_t \\
&&+\Big[ [\l^{22}_t \th^2_t + \l^{21}_t \th^1_t + l^2_t]  + \ol \G_t [\l^{12}_t \th^2_t + \l^{11}_t \th^1_t + l^1_t] + \ol \l_t \th^1_t + \ol a_t[a^{22}_t \th^2_t + a^{21}_t \th^1_t + b^2_t]^* \Big]: d\la\bo\ra_t.\nonumber
\eea
We want to choose $\ol a, \ol \l$ so that the right side above involves only $\ol \th$. That is, 
\beaa
a^{21} + \ol \G_t a^{11} + \ol a = \ol \G [a^{22}+ \ol \G a^{12}],\q \l^{21} + \ol\G \l^{11} + \ol \l + \ol a (a^{21})^* = \ol\G_t[\l^{22} + \ol \G \l^{12} + \ol a( a^{22})^*].
\eeaa
This implies
\bea
\label{ola}
\ol a &=& a^{12} (\ol \G)^2 + [a^{22}-a^{11}]\ol \G - a^{21};\\
\ol \l &=& \l^{12} (\ol \G)^2 +  [\l^{22} - \l^{11}]\ol \G -\l^{21} + \ol a [a^{22} \ol \G - a^{21}]^* \nonumber\\
&=& c^3 (\ol\G)^3 + c^2 (\ol\G)^2 + c^1 \ol\G + c^0, \q\mbox{where}\nonumber\\
c^3&:=& a^{12}(a^{22})^*,\q c^2 \;:=\; \l^{12} - a^{12} (a^{21})^* +(a^{22}-a^{11})(a^{22})^*\nonumber \\
c^1 &:=& \l^{22}-\l^{11} -(a^{22}-a^{11})(a^{21})^* - a^{21}(a^{22})^*,\q c^0\; :=\; a^{21}(a^{21})^*- \l^{21}.\nonumber
\eea
 Plugging this into \reff{barG} we obtain the following Riccati type of RDE:
\bea
\label{Riccati}
d \ol\G_t = \Big[ a^{12}_t (\ol \G)^2_t + [a^{22}_t-a^{11}_t]\ol \G_t - a^{21}_t\Big] \cd d\bo_t + \Big[c^3_t (\ol\G)^3_t + c^2_t (\ol\G)^2_t + c^1_t \ol\G_t + c^0_t\Big] : d\la\bo\ra_t,
\eea
and the RDE \reff{barth} becomes:
\bea
\label{barth2}
d\ol \th_t &=& \Big[[a^{22}+ \ol \G a^{12}] \ol\th_t + [b^2_t + \ol\G_t  b^1_t]  \Big] \cd d\bo_t \\
&&+\Big[ [\l^{22} + \ol \G \l^{12} + \ol a( a^{22})^*] \ol\th_t + [ l^2_t  + \ol \G_t l^1_t + \ol a_t (b^2_t)^*] \Big]: d\la\bo\ra_t.\nonumber
\eea
Moreover, plug $\th^2 = \ol\th - \ol\G \th^1$ into the second equation of \reff{linear2}, we have
\bea
\label{linear3}
 d \th^1_t = \Big[  [a^{11}_t - a^{12} \ol\G_t] \th^1_t  +  [a^{12}_t \ol \th_t +  b^1_t] \Big]\cd d\bo_t +\Big[ [ \l^{11}_t -\l^{12}_t\ol\G_t] \th^1_t +[\l^{12}_t \ol \th_t + l^1_t]\Big] : d\la\bo\ra_t.
\eea

Now the RDEs \reff{Riccati}, \reff{barth2}, and \reff{linear3} are decoupled. We shall emphasize though  the Riccati RDE \reff{Riccati} typically does not have a global solution on $[0, T]$.  However, following the arguments in Theorem \ref{thm-RDEwell}, there exists a constant $\d>0$, which depends only on the coefficients  $a$, $\l$ and the rough path $\bo$, such that the Riccati RDE \reff{Riccati} with initial value $0$ has a solution whenever the time interval is smaller than $\d$.  We now set $0=t_0<\cds<t_n=T$ such that $t_i - t_{i-1} \le \d$ for $i=1,\cds, n$, and we solve the system \reff{linear2} as follows. First,  we solve RDE \reff{Riccati} on $[t_0, t_1]$ with initial value $\ol\G_{t_0} = 0$.  Plug this into \reff{barth2}, where $\ol a$ is determined by \reff{ola}, we solve \reff{barth2} on $[t_0, t_1]$ with initial value $\ol\th_0 = \th^2_0$. Plug $\ol\G$ and $\ol \th$ into \reff{linear3}, we may solve \reff{linear3} on $[t_0, t_1]$ with initial value $\th^1_0$. Moreover, $\th^2 := \ol\th - \ol\G \th^1$ satisfies the second equation of \reff{linear2} on $[t_0, t_1]$ with initial value $\th^2_0$. Next, we solve the Riccati RDE \reff{Riccati} on $[t_1, t_2]$, again with initial value $\ol\G_{t_1} = 0$. Then we solve  \reff{barth2} on $[t_1, t_2]$ with initial value $\ol\th_{t_1} = \th^2_{t_1}$. Plug $\ol\G$ and $\ol \th$ into \reff{linear3}, we may solve \reff{linear3} on $[t_1, t_2]$ with initial value $\th^1_{t_1}$. Moreover, $\th^2 := \ol\th - \ol\G \th^1$ satisfies the second equation of \reff{linear2} on $[t_1, t_2]$ with initial value $\th^2_{t_1}$. Repeat the arguments we solve the system \reff{linear2} over the whole interval $[0, T]$.

{\it Step 3.} We now assume the result is true for $|E|=n-1$ and we shall prove the case $|E|=n$. With obvious notations, we consider
\bea
\label{linearn}
d \th^i_t =\Big[\sum_{j=1}^{n} a^{ij}_t \th^j_t + b^i_t\Big] \cd d\bo_t + \Big[\sum_{j=1}^{n} \l^{ij}_t \th^j_t + l^i_t\Big] : d\la\bo\ra_t,\q i=1,\cds, n.
\eea
Denote $\ol\th := \th^{n} + \sum_{i=1}^{n-1} \ol\G^i \th^i$, where, for $i=1,\cds, n-1$,
\bea
\label{olGi}
d\ol\G^i_t &=& \Big[ \sum_{j=1}^{n-1} [a^{jn} \ol\G^i_t - a^{ji}_t] \ol\G^j_t + [a^{nn}_t \ol\G^i_t - a^{ni}_t]\Big]  \cd d\bo_t\nonumber \\
&&+ \Big[[\ol \G^i_t \l^{nn}_t - \l^{ni}_t] + \sum_{j=1}^{n-1} \ol \G^j_t [\ol\G^i_t \l^{jn}_t -\l^{ji}_t]\\
&&+ \sum_{j=1}^{n-1}\big[ \sum_{k=1}^{n-1} [a^{kn} \ol\G^j_t-a^{kj}_t] \ol\G^k_t + [a^{nn}_t \ol\G^j_t  - a^{nj}_t]\big]  [\ol\G^i_t (a^{jn}_t)^*-(a^{ji}_t)^*]\Big] : d\la\bo\ra_t,\nonumber
\eea
Then
\bea
\label{linearn2}
d\ol \th_t &=& \Big[ [a^{nn}_t + \sum_{i=1}^{n-1} \ol\G^i_t a^{in}_t ]\ol\th_t + [b^n_t + \sum_{i=1}^{n-1} \ol\G^i_t b^i_t]  \Big] \cd d\bo_t \\
&&+\Big[ \big[\l^{nn}_t + \sum_{i=1}^{n-1}[\ol\G^i_t \l^{in}_t +  \ol a^i_t (a^{in})^*]\big] \ol\th_t +\big[l^n+ \sum_{i=1}^{n-1}[\ol\G^i_t l^i_t + \ol a^i_t (b^i_t)^*]\big] \Big]: d\la\bo\ra_t.\nonumber\\
\mbox{where}&& \ol a^i_t := \sum_{j=1}^{n-1} [a^{jn} \ol\G^i_t - a^{ji}_t] \ol\G^j_t + [a^{nn}_t \ol\G^i_t - a^{ni}_t].\nonumber
\eea
Plug this into \reff{linearn}, we obtain
\bea
\label{linearn3}
d \th^i_t &=& \Big[\sum_{j=1}^{n-1} [a^{ij}_t - a^{in}_t\ol\G^j_t] \th^j_t + [b^i_t + a^{in}_t \ol\th_t]\Big] \cd d\bo_t \\
&&+ \Big[\sum_{j=1}^{n-1} [\l^{ij}_t - \l^{in}_t \ol\G^j_t]\th^j_t + [l^i_t + \l^{in}_t \ol\th_t]\Big] : d\la\bo\ra_t,\q i=1,\cds, n-1.\nonumber
\eea
Now similarly, there exists $\d>0$, depending only on $a$, $\l$, and the rough path $\bo$, such that the system of Riccati type RDE \reff{olGi} with initial condition $0$ has a solution whenever the time interval is smaller than $\d$. Now set $0=t_0<\cds<t_n=T$ such that $t_i-t_{i-1}\le \d$. As in Step 2, we may first solve \reff{olGi} on $[t_0, t_1]$ with initial condition $\ol\G^i_0=0$. We then solve \reff{linearn2} on $[t_0, t_1]$ with initial condition $\ol\th_0 = \th^n_0$. Now notice that the linear system \reff{linearn3} has only dimension $n-1$, then by induction assumption, we may solve \reff{linearn3} on $[t_0, t_1]$ with initial condition $\th^i_0$, $i=1,\cds, n-1$, which further provides $\th^n := \ol\th - \sum_{i=1}^{n-1} \ol\G^i \th^i$. Now repeat the arguments as in Step 2, we obtain the solution over the whole interval $[0, T]$.
\qed

\begin{rem}
\label{rem-linear}
{\rm  (i) When $E=\dbR$, the representation formula \reff{linearrep} actually holds under weaker conditions: $a, b \in \cC^1_{\o,\ab}(\dbR^{d})$. Moreover, uniqueness also holds under this weaker condition. Indeed, for any arbitrary solution $\th\in \cC^2_{\o,\ab}(E)$ and for the $\G$ defined in \reff{linearrep},  by applying the It\^{o}-Ventzell formula \reff{Ventzell}  we see that
\beaa
\G_t  \th_t = \th_0 +  \int_0^t \G_s b_s \cd d\bo_s + \int_0^t \G_s \big[l_s - a_s b^*_s\big] : d\la \bo\ra_s.
\eeaa
Then $\th$ has to be the one in \reff{linearrep}.

(ii) In multidimensional case, we note that the Riccati RDE \reff{Riccati} does not involve $b$. Then we may also obtain the uniqueness, under our weaker condition $b\in \cC^1_{\o,\ab}(E^d)$, from the strategy in this proof.
\qed}
\end{rem}

Applying Theorem \ref{thm-RDEstability} and following the arguments in the beginning of the proof for Theorem \ref{thm-RDElinear} (or Remark \ref{rem-linear} (ii)) concerning the weaker condition on $b$, the following result is immediate.
\begin{cor}
\label{cor-RDElinear}  Let  $\bo, a, b, \l, l, \th$ be as in Theorem \ref{thm-RDElinear} and   $\tilde \bo, \tilde a, \tilde b, \tilde \l, \tilde l, \tilde \th$. Denote $\D \f := \f-\tilde \f$ for appropriate $\f$. Then
\beaa
d^{\o,\tilde\o}_\ab(\th,\tilde\th)  &\le& C\Big[d^{\o,\tilde\o}_\ab(a, \tilde a) + d^{\o,\tilde\o}_\ab(b, \tilde b) + \|\D \l\|_\b + \|\D l\|_\b + \|\D \bo\|_\a\\
&& + |\D a_0|+|\pa_\o a_0-\pa_{\tilde\o}\tilde a_0| + |\D b_0|+|\pa_\o b_0 - \pa_{\tilde\o} \tilde b_0|\Big].
\eeaa
\end{cor}

\section{Pathwise solutions of stochastic differential equations}
\label{sect-SDE}
\setcounter{equation}{0}

\subsection{The rough path setting for Brownian motion}
Let $\O_0 := \{\o \in C([0, T], \dbR^d): \o_0=0\}$ be the canonical space, $B$ the canonical process, $\dbF=\dbF^B$ the natural filtration, and $\dbP_0$ the Wiener measure. Following  F\"{o}llmer \cite{Follmer} (or see Bichteler \cite{Bichteler} and Karandikar \cite{Karandikar} for more general results on pathwise stochastic integration),  we may construct  pathwise  It\^{o} integration as follows:
\bea
\label{Phi}
\Phi_t(\o) := \limsup_{n\to \infty} \sum_{i=0}^{2^n-1} \o_{t^n_i} (\o_{t^n_i\wedge t, t^n_{i+1}\wedge t})^* &\mbox{where}&  t^n_i:= {iT\over 2^n}, i=0,\cds, 2^n.
\eea
Then $\Phi$ is $\dbF$-adapted and 
$\Phi_t = \int_0^t B_s d_{Ito} B^*_s$, $0\le t\le T$, $\dbP_0$-a.s.
Here $d_{Ito}$ stands for It\^{o} integration.
Define 
\bea
\label{uF}
\left.\ba{c}
\ul \Phi_{s,t}(\o) := \Phi_t(\o) - \Phi_s(\o) - \o_s \o_{s,t}^*,\q \ul \Phi^{Str}_{s,t}(\o) := \ul \Phi_{s,t}(\o) + {1\over 2}(t-s)I_d;\\
\la\o\ra_t := \o_t \o_t^* - \Phi_t(\o) - [\Phi_t(\o)]^*.
\ea\right.
\eea
 It is straightforward to check that
 \bea
 \label{uFproperty}
 \ul \Phi_{s,t}(\o) - \ul \Phi_{s,r}(\o) - \ul \Phi_{r, t}(\o) = \o_{s,r}\o_{r,t}^* = \ul \Phi^{Str}_{s,t}(\o) - \ul \Phi^{Str}_{s,r}(\o) - \ul \Phi^{Str}_{r, t}(\o)  .  
 \eea
 Moreover, we have the following well known result:
 \begin{lem}
 \label{lem-uF}
 For any ${1\over 3}<\a<{1\over 2}$, we have $\dbP_0( A_\a) = 1$, where
 \bea
 \label{Aa}
A_\a &:=& \Big\{\sup_{(s,t)\in \dbT^2} {| \ul \Phi_{s,t} |\over |t-s|^{2\a}} < \infty \Big\} \cap \Big\{\la\o\ra_t = t I_d, 0\le t\le T\Big\}\\
&&\cap\Big\{ \limsup_{t\downarrow s} {|v\cd \o_{s,t} |\over |t-s|^{2\a}} = \infty, \forall s\in \dbQ\cap [0, T),~ v\in \dbR^d\backslash \{0\} \Big\}.\nonumber
 \eea
 \end{lem}

 Now set, for the $A_\a$ defined in \reff{Aa},
 \bea
 \label{O}
\left.\ba{c}
  \O :=  \Big\{\o\in \O_0:   (\o, \ul \Phi(\o)) \in \bO_\a ~\mbox{and}~ \o\in A_\a,  ~\mbox{for all}~{1\over 3}<\a< {1\over 2}\Big\};\\
 d_\a(\o, \tilde\o) := d_\a\Big((\o, \ul \Phi(\o)), (\tilde\o, \ul \Phi(\tilde\o))\Big),\q \mbox{for all}~\o, \tilde\o\in \O~\mbox{and}~{1\over 3}<\a< {1\over 2}.
 \ea\right.
 \eea
 By \reff{uFproperty}  and Lemma \ref{lem-uF}, we see that $\dbP_0(\O)=1$.  From now on, we shall always restrict the sample space to $\O$, and we still denote by $B$  the canonical process and $\dbF := \dbF^B$. Define
 \bea
 \label{caO}
 &\cC(\O, E) := \bigcup \Big\{\cC_\ab(\O, E): \ab~\mbox{satisfies}~\reff{ab}\Big\},\q \mbox{where}&\\
&\cC_\ab(\O, E):=  \Big\{\th\in \dbL^0(\dbF):  \th(\o) \in  \cC^1_{\o,\ab}(E), ~\forall \o\in\O, ~\mbox{and}~ \dbE^{\dbP_0}\Big[\|\th(\o)\|_{\o,\ab}^2 \Big] <\infty\Big\}.&\nonumber
\eea
We now define pathwise stochastic integral by using rough path integral:  for $\th \in \cC(\O, E^d)$, 
\bea
\label{pI}
\left.\ba{c}
\dis\Big(\int_0^t \th_s \cd dB_s\Big)(\o) := \int_0^t \th_s(\o) \cd d (\o, \ul \Phi(\o))_s,\q\forall \o\in \O;\\
\dis\Big(\int_0^t \th_s \circ dB_s\Big)(\o) := \int_0^t \th_s(\o) \cd d (\o, \ul \Phi^{Str}(\o))_s,\q\forall \o\in \O.
\ea\right.
\eea
The following result can be found in \cite{FH} Proposition 5.1 and Corollary 5.2.
\begin{thm}
\label{thm-pI}
 For any $\th \in \cC(\O, E^d)$, the above pathwise stochastic integrals $\int_0^t \th_s \cd dB_s$ and $\int_0^t \th_s \circ dB_s$ coincide with the It\^{o} integral and the Stratonovic integral, respecively.
 \end{thm}

\begin{rem}
\label{rem-pI}
{\rm 
Let $X$ be a semi-martingale with $d X_t = \th_t \cd dB_t + \l_t dt$, where $\th\in \cC(\O, E^d)$ and $\l$ is continuous.  Then $X\in \cC(\O, E)$ with $\pa_\o X_t(\o) = \th_t(\o)$ for each $\o\in \O$. In the spirit of Dupire \cite{Dupire}'s functional It\^{o} calculus, \cite{BMZ} defines the above $\th$ as the path derivative of the process $X$. So the Gubinelli's derivative $\pa_\o X(\o)$ in Definition \ref{defn-Gubinelli}  is consistent with the path dervatives introduced in \cite{BMZ}.
\qed}
\end{rem}

\begin{rem}
\label{rem-DoBM}
{\rm Let $\o\in \O$ and $\th \in \cC^2_{(\o, \ul \Phi(\o)), \ab}(E)$ for certain $\ab$ satisfying \reff{ab}. Define
\bea
\label{paot}
\pa^\o_t \th := \mbox{Trace} (D^\o_t \th).
\eea
Then $\pa^\o_t \th$ is unique and is consistent with the time derivative in \cite{BMZ}. Moreover, the pathwise Ito formula \reff{Ito} and the pathwise Taylor expansion \reff{Taylor1}, \reff{Taylor2} become:
\bea
\label{ItoTaylor}
d \th_t \!\!\! &=& \!\!\! \pa_\o \th_t d\bo_t + \Big[\pa^\o_t \th_t + {1\over 2} \mbox{Trace}(\pa^2_{\o\o}\th_t)\Big]dt;\nonumber\\
 \th_{s,t} \!\!\! &=&\!\!\! \pa_\o \th_s \o_{s,t} + {1\over 2} \pa^2_{\o\o} \th_s  : [\o_{s,t} \o^*_{s,t} + \uo_{s,t} - \uo_{s,t}^*] + \pa^\bo_t \th_s (t-s)  + O((t-s)^{2\a+\b});\\
 \th_{s,t} \!\!\! &=&\!\!\! \pa_\o \th_s \o_{s,t} + {1\over 2} \pa^2_{\o\o} \th_s  : [\o_{s,t} \o^*_{s,t}] + \pa^\bo_t \th_s (t-s)  + O((t-s)^{2\a+\b}),\nonumber
 \eea
respectively. These are also consistent with \cite{BMZ}.
\qed}
\end{rem}

\subsection{Stochastic differential equations with regular solutions} 

We now consider the following SDE  with random coefficients:
\bea
\label{SDE}
X_t = x  + \int_0^t \si(s, X_s,\o)\cd  dB_s + \int_0^t b(s, X_s,\o) ds,\q\o\in \O,
\eea
where $b, \si$ are $\dbF$-progressively measurable. Clearly, the above SDE can be rewritten as the following RDE:
\bea
\label{SDE-RDE}
X_t (\o) = x + \int_0^t \si(s, X_s(\o),\o) \cd d(\o, \ul \Phi(\o))_s + \int_0^t b(s, X_s(\o),\o) {I_d\over d} : d\la \o\ra_s ,~\o\in \O.
\eea
The following result is a direct consequence of Theorems \ref{thm-RDEwell} and \ref{thm-RDEstability}.

\begin{thm}
\label{thm-SDE} (i) Assume, for each $\o\in \O$,  there exists $\ab(\o)$ satisfying \reff{ab} such that $b(\cd, \o)\in \cC^2_{\b(\o)}(E, E)$ and $\si(\cd, \o) \in \cC^{2,3}_{\o,\ab(\o)}(E, E^d)$. Then the SDE has a unique solution $X$ such that $X(\o)\in \cC^2_{\o,\ab(\o)}(E)$ for all $\o\in\O$.

(ii) Assume further that $b$ and $\si$ are continuous in $\o$ in the following sense: 
\bea
\label{bsicont}
&\lim_{n\to \infty} \Big[\|b(\cd, \o^n) - b(\cd, \o)\|_{1, \b(\o)} + d_{2,\ab(\o)}^{\o, \o^n}(\si(\cd, \o^n),  \si(\cd, \o)) \Big] = 0,&\\
&\mbox{for any $\o, \o^n\in \O$ such that} ~ \lim_{n\to \infty} d_{\a(\o)}(\o^n, \o) = 0.&\nonumber
\eea
Then $X$ is also continuous in $\o$ in the sense that:
\bea
\label{Xcont}
\lim_{n\to \infty} d_{\ab(\o)}^{\o, \o^n}(X(\o), X(\o^n)) = 0, ~~\mbox{and consequently,}~~
\lim_{n\to \infty}  \|X(\o)- X(\o^n)\|_\infty = 0.
\eea
\end{thm}

\bs

\begin{rem}
\label{rem-SDE0}
{\rm  The construction of pathwise solutions of SDEs via rough path is standard. However, we remark that our canonical sample space $\O$ is universal, which particularly does not depend on the controlled rough path $\th$ or the coefficient $\si(t,\o,x)$. Consequently, our solution is constructed indeed for every $\o\in \O$, without the exceptional null set.  To our best knowledge, such a message is new.
\qed}
\end{rem}

\begin{rem}
\label{rem-SDE}
{\rm (i) Assume $\si$ is H\"{o}lder-${1\over 2}$ continuous in $t$ and Lipschitz continuous in $\o$ in the following sense:
\bea
\label{siLip}
|\si(t, x, \o) - \si(\tilde t, x, \tilde\o)| \le C\Big[\sqrt{\tilde t - t} + \sup_{0\le s\le T}|\o_{s\wedge t} - \tilde \o_{s\wedge \tilde t}|\Big],
\eea
Then $\si(\cd, x, \o)$ is H\"{o}lder-$\a$ continuous in $t$ for all $\a<{1\over 2}$. We remark that the distance in the right side of \reff{siLip} is used in Zhang and Zhuo \cite{ZZ} and is equivalent to the metric introduced by Dupire \cite{Dupire}.

(ii) As mentioned in Introduction, since $\o$ is only H\"{o}lder-$\a$ continuous for $\a<{1\over 2}$, it is not reasonable to assume $\si(\cd, x, \o)$ is H\"{o}lder-$(1-\a)$ continuous as required in Lejay and Victoir \cite{LV}. 
\qed}
\end{rem}

\begin{rem}
\label{rem-Stratonovich}
{\rm Under the Stratonovich integration, the quadratic variation of Brownian motion sample path vanishes: $\la (\o, \ul \Phi^{str}(\o))\ra_t = 0$. If we want to consider SDE in the form:
\bea
\label{SDE-Stratonovich}
d X_t =\si(t,X_t,\o) \circ dB_t + b(t, \o, X_t) dt,
\eea
we cannot simply rewrite it into
\beaa
d X_t(\o) = \si(t,\o,X_t(\o)) \cd d(\o, \ul \Phi^{str}(\o))_t + b(t, \o, X_t(\o)) {I_d\over d} : d\la  (\o, \ul \Phi^{str}(\o))\ra_t.
\eeaa
 We can obtain pathwise solution of \reff{SDE-Stratonovich} in the following two ways:

(i) We may rewrite \reff{SDE-Stratonovich} in It\^{o} form:
\bea
\label{Stratonovich-Ito}
d X_t = \si(t,\o,X_t) \cd dB_t +  \Big[b + {1\over 2}\mbox{Trace}\big(\pa_\o\si + \pa_y\si \otimes \si^*\big) \Big](t, \o, X_t) dt,
\eea
which corresponds further to the following RDE:
\bea
\label{Stratonovich-RDE}
d X_t(\o) \!= \! \si(t,\o,X_t(\o)) \cd d(\o, \ul\Phi(\o))_t +  \Big[{b I_d\over d} + {\pa_\o\si + \pa_y\si \otimes \si^*\over 2}\Big](t, \o, X_t(\o)) : d\la\o\ra_t.
\eea

(ii) In Section \ref{sect-RDE}, we may easily extend our results to more general RDEs:
\bea
\label{RDEdt}
d\th_t = g(t,\th_t) \cd d\bo_t + f(t,\th_t) : d\la\bo\ra_t + h(t,\th_t) dt.
\eea
Then we may deal with \reff{SDE-Stratonovich} directly.
\qed}
\end{rem}

\section{Rough PDEs and Stochastic PDEs}
\label{sect-RPDE}
\setcounter{equation}{0}

In this section, we extend the results in previous sections to rough PDEs \reff{RPDE0} and stochastic PDEs \reff{SPDE0} with random coefficients. The wellposedness of such RPDEs and SPDEs, especially in fully nonlinear case, is very challenging and has received very strong attention. We refer to Lions and Souganidis \cite{LS1, LS2, LS3, LS4}, Buckdahn and Ma \cite{BM1, BM2}, Caruana and Friz \cite{CF1}, Caruana,  Friz and Oberhauser \cite{CFO}, Friz and Obhauser \cite{FO},   Diehl and  Friz \cite{DF},  Oberhauser and Riedel \cite{DOR}, and Gubinelli,  Tindel and Torrecilla,  \cite{GTT} for wellposedness of some classes of RPDEs/SPDEs, where  various notions of solutions are proposed.

While this section is mainly motivated from the study of pathwise viscosity solutions of SPDEs in Buckdahn, Ma  and Zhang \cite{BMZ2},  in this section we shall focus on calssical solutions only. In particular, we do not intend to establish strong wellposedness for general $f$, instead we shall investigate diffusion coefficients $\si$ and $g$ and see when the RPDE/SPDE can be transformed to a deterministic PDE. Again, unlike most results in standard literature of rough PDEs, we allow the coefficients to depend on $(t, \o)$. The results will require quite high regularity of the coefficients, in the sense of our path regularity. In order to simplify the presentation, for some results we shall not specify the precise regularity conditions.

\subsection{RDEs with spatial parameters}
\label{sect-RDEx}
Let $u_0: \tilde E\to E$, $g: \dbT \times \tilde E\times  E\to E^d$, $f: \dbT \times \tilde E \times E\to E^{d\times d}$, and consider the following RDE with parameter $x\in \tilde E$:
\bea
\label{RDEu}
u_t(x) = u_0(x) + \int_0^t g(s, x, u_s(x)) \cd d\bo_s + \int_0^t f(s, x, u_s(x)): d\la\bo\ra_s.
\eea
Assume  $u_0$, $g$ and $f$ are differentiable in $x$,  and differentiate \reff{RDEu} formally in $x_i$, $i=1,\cds, |\tilde E|$, we obtain: denoting  $v^i_t(x) := \pa_{x_i} u_t(x)$, 
\bea
\label{RDEv}
&\dis v^i_t(x) = \pa_{x_i} u_0(x) + \int_0^t [\pa_{x_i} g (s, x, u_s(x)) +  \pa_y g(s, x, u_s(x)) \otimes v^i_s(x)] \cd d \bo_s&\nonumber\\
&\dis\qq\qq\qq + \int_0^t [\pa_{x_i} f(s, x, u_s(x)) +  \pa_y f(s, x, u_s(x))\otimes v^i_s(x)] : d\la\bo\ra_s.&
\eea

\begin{thm}
\label{thm-thx}
 Assume

(i) $u_0$, $g$, $f$ are continuously differentiable in $x$;

(ii)  for each $x\in \tilde E$, $i=1,\cds, |\tilde E|$, $j=1,\cds, |E|$,
\bea
\label{thx}
\left.\ba{c}
g(x,\cd) \in \cC^{2,3}_{\o,\ab}(E, E^d),\q    f(x,\cd) \in \cC^2_\b(E, E^{d\times d});\\
\pa_{x_i} g(x,\cd) \in \cC^{1,2}_{\o,\ab}(E, E^d),\q  \pa_{y_j} g(x,\cd) \in \cC^{2,3}_{\o,\ab}(E, E^d),\q  \pa_{x_i} f(x,\cd)\in \cC^0_\b(E, E^{d\times d}).
\ea\right.
\eea

(iii) 
 for any $x\in \tilde E$, denoting $\D \f := \f(x+\D x,\cd) - \f(x,\cd)$ for appropriate $\f$,
\bea
\label{hcont2}
\left.\ba{c}
\lim_{|\D x|\to 0} \Big[\|\D g\|_{2,\o,\ab} + \|\D  f\|_{1,\b}\Big] = 0;\\
\lim_{|\D x|\to 0} \Big[ \|\D[\pa_x g]\|_{2,\o,\ab} + \|\D[\pa_y g]\|_{2,\o,\ab}+ \|\D [\pa_x f]\|_{0,\b} + \|\D [\pa_y f]\|_{0,\b}\Big] = 0.
\ea\right.
\eea
Moreover, $\pa_{\o x} g$ and $\pa_{\o y}g$ are continuous.

 Then, for each $x\in \tilde E$, RDEs \reff{RDEu} and \reff{RDEv} have unique solution $u(x,\cd), v^i(x,\cd) \in \cC^2_{\o,\ab}(E)$, respectively. Moreover, $u$ is differentiable  in $x$ with $\pa_{x_i} u = v^i$. 
\end{thm}
\proof First, without loss of generality we may assume $|\tilde E|=1$, namely $\tilde E= \dbR$. For each $x\in \tilde E$, by the first line of \reff{thx} and  applying Theorem \ref{thm-RDEwell}, we see that RDE \reff{RDE} has a unique solution $u(x) \in \cC^2_{\o,\ab}(E)$.  By the second line of \reff{thx} and applying Theorem \ref{thm-chain} and Lemma \ref{lem-Youngf}, we see that, for $j=1,\cds, |E|$,
\beaa
\pa_{x} g (x, u(x))\in \cC^1_{\o,\ab}(E^d), \pa_{y_j}g (x, u(x)) \in \cC^2_{\o,\ab}(E^d),   \pa_{x} f (x, u(x)), \pa_{y_j} f(x, u(x))\in \O_\b(E^{d\times d}).
\eeaa
Then by  Theorem  \ref{thm-RDElinear} the linear RDE \reff{RDEv} has a unique solution $v(x) \in \cC^2_{\o,\ab}(E)$. 

It remains to prove $\pa_x u = v$. Given $x\in \dbR$, $\D x\in \dbR\backslash \{0\}$ and $\l \in [0,1]$, denote
\beaa
&\D u_t := u_t(x+\D x) - u_t(x),\q \td u_t := {\D u_t \over \D x},&\\
&\f_t(\l) := \f(t, x+\l \D x, u_t(x) + \l \D u_t(x)),~ \D \f_t(\l) := \f_t(\l) - \f_t(0),~\mbox{for appropriate $\f$}.&
\eeaa
By the first line of \reff{hcont2}, it follows from Theorem \ref{thm-RDEstability} that: 
\bea
\label{Duconv}
\lim_{|\D x|\to 0} \|\D u\|_{\o,\ab} = 0.
\eea
Moreover, one can easily check that,
\beaa
d\td u_t &=& \int_0^1 [\pa_x g_t(\l) + \pa_y g_t(\l) \otimes \td u_t] d\l \cd d\bo_t +\int_0^1 [\pa_x f_t(\l) + \pa_y f_t(\l) \otimes \td u_t] d\l : d\la\bo\ra_t;\\
 d v_t(x) &=& [\pa_x g_t(0) + \pa_y g_t(0) \otimes v_t(x)] \cd d\bo_t +[\pa_x f_t(0) + \pa_y f_t(0) \otimes v_t(x)] : d\la\bo\ra_t.
\eeaa 
By the  second line of \reff{hcont2} and  \reff{Duconv}, it follows from Lemmas \ref{lem-g} (ii) and  \ref{lem-Youngf} (ii) that
\beaa
&\lim_{|\D x|\to 0} \Big[ \|  \pa_x g_t(\l) - \pa_x g(0)\|_{\o,\ab} +  \|  \pa_y g_t(\l) - \pa_y g(0)\|_{\o,\ab} \Big] = 0;&\\
&\lim_{|\D x|\to 0} \Big[  \|  \pa_x f_t(\l) - \pa_x f(0)\|_{\b} +\|  \pa_y f_t(\l)- \pa_y f(0)\|_{\b}\Big] = 0.&
\eeaa
 for any $\l\in [0,1]$. Furthermore, by Theorem \reff{thm-chain} (i) we have
 \beaa
 \pa_\o [\pa_x g_0(\l)] = \pa_{\o x} g(\l) + \pa_{yx} g_0(\l) \otimes g_0(\l),\q \pa_\o [\pa_y g_0(\l)] = \pa_{\o y} g(\l) + \pa_{yy} g_0(\l) \otimes g_0(\l)
 \eeaa
 Recalling the continuity of $\pa_{\o x} g$, $\pa_{\o, y} g$ in (iii) we see that,   for any $\l\in [0,1]$,
 \beaa
 \lim_{|\D x|\to 0} \Big[ |  \pa_\o [\pa_x g_0(\l)]  -\pa_\o [\pa_x g_0(\l)| +  |  \pa_\o [\pa_y g_0(\l)]  -\pa_\o [\pa_y g_0(\l)| \Big] = 0.
 \eeaa
 Now by Corollary \ref{cor-RDElinear} we have $\lim_{|\D x|\to 0}\|\td u - v(x)\|_{\o,\ab}=0$.
That is, $\pa_x u_t(x) = v_t(x)$.
 \qed

\subsection{Pathwise characteristics}
As standard in the literature, see e.g. Kunita \cite{Kunita} for Stochastic PDEs and \cite{FH} Chapter 12 for rough PDEs, the main tool for dealing with semilinear RPDEs/SPDEs is the characteristics, which we shall introduce below by using RDEs against rough paths and backward rough paths. 

Let $\si: \dbT \times \tilde E \to \tilde E^d$ and $g: \dbT \times \tilde E \times E \to E^{d\times d}$. Fix $t_0\in \dbT$ and denote
\bea
\label{Backsi}
\lsi^{t_0}(t, y) := \si(t_0-t, y), \q \Lg^{t_0}(t, x, y) := g(t_0-t, x, y),
\eea
Consider  the following characteristic RDEs:
\bea
\label{BackRDE1}
 &\dis\th^x_t =x - \int_0^t \si(s, \th^x_s) \cd d \bo_s,\q \lth^{t_0,x}_t = x + \int_0^t \lsi^{t_0}(s, \lth^{t_0,x}_s) \cd d\lbo^{t_0}_s;&\\
\label{BackRDE2}
&\dis \eta^{x,y}_t =y + \int_0^t g(s, \th^x_s, \eta^{x,y}_s) \cd d \bo_s,\q \leta^{t_0,x,y}_t = y - \int_0^t \Lg^{t_0}(s,  \lth^{t_0,x}_s, \leta^{t_0,x,y}_s) \cd d\lbo^{t_0}_s.&
\eea
By Lemma \ref{lem-Back} and Theorem \ref{thm-RDEwell}, the following result is obvious.
\begin{lem}
\label{lem-BackRDE} (i) Assume $\si \in \cC^{2,3}_{\o,\a}(\tilde E, \tilde E^d)$. Then, for each $x\in \tilde E$,  the RDEs \reff{BackRDE1} have unique solution $\th^x\in \cC^1_{\o, \ab}(\tilde E)$ and $\lth^{t_0,x}\in \cC^1_{\lo^{t_0}, \ab}([0, t_0], \tilde E)$ satisfying $\lth^{t_0, \th^x_{t_0}}_t = \th^x_{t_0-t}$, $t\in [0, t_0]$. In particular, the mapping $x \mapsto \th^x_{t_0}$  is one to one with inverse function $x \mapsto \lth^{t_0,x}_{t_0}$.

(ii) Assume further that, for each $x\in \tilde E$ and for the above solution $\th^x$,  the mapping $(t, y) \mapsto g(t, \th^x_t, y)$ is in $ \cC^{2,3}_{\o,\a}(E, E^{d\times d})$. Then the RDEs \reff{BackRDE2} have unique solution $\eta^{x,y} \in  \cC^1_{\o,\ab}(E)$ and $\leta^{t_0,x,y} \in \cC^1_{\lo^{t_0},\ab}(E)$ satisfying $\leta^{t_0, \th^x_{t_0}, \eta^x_{t_0}}_t = \eta^{x,y}_{t_0-t}$, $t\in [0, t_0]$. In particular, the mapping $(x, y) \mapsto (\th^x_{t_0}, \eta^{x,y}_{t_0})$  is one to one with inverse functions $(x, y) \mapsto (\lth^{t_0,x}_{t_0}, \leta^{t_0,x,y}_{t_0})$.
\end{lem}

Now define
\bea
\label{f}
\f(t,x) := \lth^{t,x}_t,\q \psi(t,x,y) :=  \leta^{t,\th^x_t,y}_t,\q \zeta(t,x,y) :=  \eta^{\f(t,x), y}_t,\q \wh g(t,x,y) := g(t,\th^x_t, y).
\eea
\begin{lem}
\label{fdiff}
Assume $\si$ and $g$ are smooth enough in the sense of Theorem \ref{thm-thx}. Then $\f, \psi$ are twice differentiable in $(x, y)$, and for any fixed $(x,y)$, $\f(\cd, x), \psi(\cd, x,y) \in \cC^\o_\a$. Moreover, they satisfy the following RDEs: 
\beaa
\f(t,x) &=& x + \int_0^t \pa_x \f\otimes \si(s,x) \cd d \bo_s \\
&&+\int_0^t \Big[{1\over 2} \pa^2_{xx} \f \otimes_2 [\si,\si] + \pa_x \f  \otimes [\pa_x \si \otimes \si^*]\Big](s,x) : d\la\bo\ra_s;\\
 \psi(t,x,y) &=& y - \int_0^t [\pa_y \psi\otimes \wh g](s, x, y) \cd d \bo_s\\
 &&+\int_0^t\Big[{1\over 2} \pa^2_{yy}\psi \otimes_2 [\wh g, \wh g] + \pa_y \psi\otimes [\pa_y \wh g \otimes \wh g^*]  \Big](s,x,y): d\la\bo\ra_s.
\eeaa
\end{lem}
\proof By  Theorem \ref{thm-thx},  $\th^x$, $\lth^{t,x}$, $\eta^{x,y}$, $\leta^{t,x,y}$ are sufficiently differentiable in $(x, y)$. This implies the desired differentiability of  $\f, \psi$. We now check the RDEs.

First, given $(s,t)\in \dbT^2$ and denote $\d := t-s$.  Note that 
\beaa
\f(t,x)=\lth^{t,x}_t = \lth^{s, \lth^{t,x}_\d}_s = \f(s, \lth^{t,x}_\d);
\eeaa
and that, applying Lemma \ref{lem-Back},
\beaa
\lth^{t,x}_\d-x&=& \int_0^\d \lsi^t(r, \lth^{t,x}_r) \cd d \lbo^t_r \\
&=& \lsi^t(0, x) \cd \lo^t_{0,\d} + [\pa_{\lo^t} \lsi^t + \pa_x\lsi^t\otimes (\lsi^t)^*](0,x): \luo^t_{0,\d} + O(\d^{2\a+\b})\\
&=& \si(t,x) \cd \o_{s,t} +[-\pa_\o \si+ \pa_x\si \otimes \si^*](t,x) :[\o_{s,t}\o^*_{s,t} - \uo_{s,t} ] + O(\d^{2\a+\b})\\
&=&\si(s,x) \cd \o_{s,t} +\pa_\o\si(s,x) :\uo_{s,t} + \pa_x\si \otimes \si^*(s,x) :[\o_{s,t}\o^*_{s,t} - \uo_{s,t} ] + O(\d^{2\a+\b})
\eeaa
Then, applying Taylor expansion,
\beaa
&&\f(t,x) -\f(s, x) = \f(s,  \lth^{t,x}_\d) - \f(s,x) \\
&=& \pa_x \f(s,x) \otimes [\lth^{t,x}_\d - x] + {1\over 2} \pa^2_{xx}\f(s,x) \otimes_2 [\lth^{t,x}_\d - x, \lth^{t,x}_\d - x ] + O(\d^{3\a})\\
&=& \pa_x \f(s,x) \otimes \Big[ \si(s,x) \cd \o_{s,t} +\pa_\o \si (s,x): \uo_{s,t} + \pa_x\si \otimes \si^*(s,x) :[\o_{s,t}\o^*_{s,t}-\uo_{s,t} ] \Big] \\
&&+ {1\over 2} \pa^2_{xx}\f(s,x) \otimes_2 [ \si(s,x) \cd \o_{s,t}] + O(\d^{2\a+\b})
\eeaa
In particular, this implies 
\beaa
\pa_\o \f = \pa_x \f \otimes \si.
\eeaa
On the other hand, by applying Theorem \ref{thm-thx} on \reff{BackRDE1} and view $(\th^x, \pa_x \th^x)$ as the solution to a higher dimensional RDE, one can check similarly that
\beaa
\pa_\o [\pa_x \f] = \pa_x[(\pa_x\f \otimes \si)^*].
\eeaa

Denote $\tilde \f$ as the right side of the RDE for $\f$. Then, taking values at $(s,x)$,
\beaa
[\tilde\f(\cd,x)]_{s,t}&=& \pa_x \f\otimes \si \cd \o_{s,t} + \pa_\o [\pa_x\f \otimes \si] : \uo_{s,t} \\
&&+ \Big[{1\over 2} \pa^2_{xx} \f \otimes [\si,\si] + \pa_x \f  \otimes [\pa_x \si \otimes \si^*]\Big] : \la\bo\ra_{s,t} + O(\d^{2\a+\b})\\
&=& \pa_x \f\otimes \si \cd \o_{s,t} + \Big[\big[\pa_x[\pa_x\f \otimes \si] \otimes \si^* + \pa_x\f \otimes \pa_\o\si\big] : \uo_{s,t} \\
&&+ \Big[{1\over 2} \pa^2_{xx} \f \otimes \si + \pa_x \f  \otimes [\pa_x \si \otimes \si^*]\Big]  : [\o_{s,t}\o^*_{s,t} - \uo_{s,t}-\uo_{s,t}^*]+ O(\d^{2\a+\b}).
\eeaa
It is straightforward to check that
$
[\f(\cd,x)]_{s,t} = [\tilde\f(\cd,x)]_{s,t} + O(\d^{2\a+\b}),
$
impling $\f = \tilde\f$.

Similarly, notice that 
\beaa
\psi(t, x, y) = \leta^{t,\th^x_t,y}_t = \leta^{s, \lth^{t,\th^x_t}_\d, \leta^{t,\th^x_t,y}_\d}_s = \leta^{s, \th^x_s, \leta^{t,\th^x_t,y}_\d}_s = \psi(s, x,  \leta^{t,\th^x_t,y}_\d).
\eeaa
Following similar arguments one can verify the RDE for $\psi$. 
\qed

\subsection{Rough PDEs}
Now consider RPDE:
\bea
\label{RPDE}
 u_t(x) &=& u_0(x) + \int_0^t [\pa_x u_s(x)\otimes \si_s(x) + g_s(x, u_s(x))] \cd d\bo_s \\
 && + \int_0^t f_s(x, u_s(x), \pa_x u_s(x), \pa^2_{xx} u_s(x)) : d\la\bo\ra_s.\nonumber
\eea
Define
\beaa
v(t,x) :=  \psi(t, x, u(t,\th^x_t))  ~\mbox{and equivalently} ~ u(t,x) = \zeta(t, x, v(t,\f(t,x))).
\eeaa
\begin{thm}
\label{thm-RPDE}
Assume the coefficients and $u$ are smooth enough. Then $u$ is a solution of RPDE \reff{RPDE} if and only if $v$ satisfies:
\bea
\label{RPDEv}
&dv_t(x) =\wh f(t,x, v_t(x), \pa_x v_t(x), \pa^2_{xx} v_t(x)):d\la \bo\ra_t,&\\
&\mbox{or equivalently}, ~D^\o_t v_t(x) = \wh f(t,x, v_t(x), \pa_x v_t(x), \pa^2_{xx} v_t(x)),&\nonumber
\eea
where
\bea
\label{RPDEF}
\wh f(t,x,y,z,\g) &:=& \pa_y \psi(t,x, \wh y) \Big[ f(t, \th^x_t, \wh y, \wh z, \wh \g)  - {1\over 2} \wh \g : [\si, \si](t,\th^x_t)\nonumber\\
&&- \big[\wh z\otimes \pa_x \si + \pa_x g + \pa_y g\otimes \wh z]\otimes \si^*\Big](t,\th^x_t, \wh y);\\
\wh y &=& \zeta(t,\th^x_t,y);\nonumber\\
\wh z &=&  \pa_x \zeta (t,\th^x_t,y) + \pa_y \zeta (t,\th^x_t,y) \otimes z \otimes \pa_x \f (t,\th^x_t);\nonumber\\
\wh \g&=& \pa^2_{xx}\xi (t,\th^x_t,y)+ [\pa_{xy}\zeta(t,\th^x_t,y)+\pa_{yx}\si(t,\th^x_t)] \otimes_2 [z, \pa_x \f(t,\th^x_t)] \nonumber\\
&&+ \pa^2_{yy}\zeta (t,\th^x_t,y) \otimes_2 [\pa_x \f\otimes \pa_x \f, \pa_x \f\otimes \pa_x\f](t,\th^x_t)\nonumber \\
&&+ \pa_y \zeta(t,\th^x_t,y) \otimes \Big[\g\otimes_2[\pa_x \f, \pa_x\f](t,\th^x_t)+ z \otimes \pa^2_{xx}\f(t,\th^x_t)\Big].\nonumber
\eea
\end{thm}
\proof  
 Applying the It\^{o}-Ventzell formula \reff{Ventzell} we have
\bea
\label{dv}
&&d u(t, \th^x_t) = g(t, \th^x_t, u(t, \th^x_t)) d\bo_t + \Big[f(\cd, u, \pa_x u, \pa^2_{xx} u) \nonumber\\
&&\q -[{1\over 2}\pa^2_x u: [\si,\si] + \pa_x u \otimes \pa_x \si\otimes \si^* + \pa_x g(\cd, u)\otimes \si^*+ \pa_y g \otimes \pa_x u\otimes \si^*]\Big](t, \th^x_t) : d\la\bo\ra_t;\nonumber\\
&&d v(t,x) = d [\psi(t,x, u(t, \th^x_t))] =\pa_y \psi(t,x, u(t, \th^x_t)) \Big[ f(\cd, u, \pa_x u, \pa^2_{xx} u)  \\
&&\qq - {1\over 2} \pa^2_x u : [\si, \si]- \big[\pa_x u\otimes \pa_x \si + \pa_x g + \pa_y g\otimes \pa_x u]\otimes \si^*\Big](t,\th^x_t, u(t,\th^x_t)):d\la\bo\ra_t.\nonumber
\eea
Now note that
\beaa
u(t,x) &=& \zeta(t,x,v(t, \f(t,x)));\\
\pa_x u &=&  \pa_x \zeta + \pa_y \zeta \otimes \pa_x v \otimes \pa_x \f;\\
\pa^2_{xx} u &=& \pa^2_{xx}\xi + [\pa_{xy}\xi+\pa_{yx}\si] \otimes_2 [\pa_x v, \pa_x \f] + \pa^2_{yy}\zeta \otimes_2 [\pa_x \f\otimes \pa_x \f, \pa_x \f\otimes \pa_x\f] \\
&&+ \pa_y \zeta \otimes \pa^2_{xx} v \otimes_2[\pa_x \f, \pa_x\f] + \pa_y\zeta \otimes\pa_x v \otimes \pa^2_{xx}\f.
\eeaa
Then
\beaa
u(t,\th^x_t) &=& \zeta(t,\th^x_t,v(t, x));\\
\pa_x u(t, \th^x_t) &=&  \pa_x \zeta (t,\th^x_t,v(t, x)) + \pa_y \zeta (t,\th^x_t,v(t, x)) \otimes \pa_x v(t,x) \otimes \pa_x \f (t,\th^x_t);\\
\pa^2_{xx} u (t,\th^x_t) &=& \pa^2_{xx}\xi (t,\th^x_t,v(t, x))+ [\pa_{xy}\zeta(t,\th^x_t,v(t, x))\\
&&+\pa_{yx}\si(t,\th^x_t)] \otimes_2 [\pa_x v(t,x), \pa_x \f(t,\th^x_t)] \\
&&+ \pa^2_{yy}\zeta (t,\th^x_t,v(t, x)) \otimes_2 [\pa_x \f\otimes \pa_x \f, \pa_x \f\otimes \pa_x\f](t,\th^x_t) \\
&&+ \pa_y \zeta(t,\th^x_t,v(t, x)) \otimes \pa^2_{xx} v (t,x)\otimes_2[\pa_x \f, \pa_x\f](t,\th^x_t) \\
&&+ \pa_y \zeta(t,\th^x_t,v(t, x)) \otimes \pa_x v(t,x) \otimes \pa^2_{xx}\f(t,\th^x_t).
\eeaa
Plug this into \reff{dv}, we obtain the result immediately.
\qed

\subsection{Pathwise solution of Stochastic PDEs}
We now study Stochastic PDE:
\bea
\label{SPDE}
 u_t(\o, x) &=& u_0(x) + \int_0^t [\si_s(\o, x) \pa_x u_s(\o,x) + g_s(\o,x, u_s(\o, x))] \cd dB_s \\
 && + \int_0^t f_s(\o, x, u_s(\o, x), \pa_x u_s(\o, x), \pa^2_{xx} u_s(\o, x))  ds,\q \dbP_0\mbox{-a.s.}\nonumber
\eea
Clearly, this corresponds to RPDE: 
\bea
\label{SRPDE}
 u_t(\o, x) &=& u_0(x) + \int_0^t [\si_s(\o, x) \pa_x u_s(\o, x) + g_s(\o, x, u_s(\o, x))] \cd d(\o, \ul F(\o))_s \\
 && + \int_0^t F_s(\o, x, u_s(\o, x), \pa_x u_s(\o,x), \pa^2_{xx} u_s(\o, x))  : d\la\o\ra_s,\q \forall \o\in \O, \nonumber\\
\mbox{where}&& F(t,\o, x,y,z,\g) := f(t,\o,x,y,z,\g){I_d\over d}.
\eea
Define $\th^{\o,x}_t$, $\psi(t,\o,x,y)$, $\wh F(t,\o,x,y,z,\g)$ in obvious sense and
\bea
\label{vto}
v(t, \o,x) := \psi(t,\o, x, u(t, \o, \th^{\o,x}_t)),\q \wh f(t,\o, x,y,z,\g) := \mbox{Trace}[\wh F(t,\o,x,y,z,\g)].
\eea
 Then we have, recalling $\pa^\o_t v$ defined in Remark \ref{rem-DoBM},
 \beaa
 d v(t,\o,x) = \pa^\o_t v(t,\o, x) dt =  \wh f_t(\o, x, v_t(\o,x), \pa_x v_t(\o,x), \pa_{xx}^2 v_t(\o,x)) dt.
 \eeaa
Clearly, this implies that $\pa^\o_t v_t(x) = \pa_t v(t,\o, x)$, the standard time derivative for fixed $(\o,x)$.  We now conclude the paper with the following result:
\begin{thm}
\label{thm-SPDE} 
Assume the coefficients and $u$ are smooth enough. Then, for each $\o\in \O$, $u(\o,\cd)$ is a  solution of \reff{SRPDE} if and only if $v(\o,\cd)$ is a solution of the following PDE:
\bea
\label{PDE} 
\pa_t v_t(\o,x) = \wh f_t(\o, x, v_t(\o,x), \pa_x v_t(\o,x), \pa_{xx}^2 v_t(\o,x)).
\eea 
\end{thm}

\end{document}